# Finding Approximate Local Minima
# Faster than Gradient Descent

Naman Agarwal
namana@cs.princeton.edu
Princeton University

Zeyuan Allen-Zhu
zeyuan@csail.mit.edu
Institute for Advanced Study

Brian Bullins
bbullins@cs.princeton.edu
Princeton University

Elad Hazan
ehazan@cs.princeton.edu
Princeton University

Tengyu Ma
tengyu@cs.princeton.edu
Princeton University

November 3, 2016

## Abstract

We design a non-convex second-order optimization algorithm that is guaranteed to return an *approximate* local minimum in time which scales linearly in the underlying dimension and the number of training examples. The time complexity of our algorithm to find an approximate local minimum is even faster than that of gradient descent to find a critical point. Our algorithm applies to a general class of optimization problems including training a neural network and other non-convex objectives arising in machine learning.

# 1 Introduction

Finding a global minimizer of a non-convex optimization problem is NP-hard. Thus, the standard goal of efficient non-convex optimization algorithms is instead to find a local minimum. This problem has become increasingly important as the state-of-the-art in machine learning is attained by non-convex models, many of which are variants of deep neural networks. Experiments in [10, 11, 21] suggest that fast convergence to a local minimum is *sufficient* for training neural nets, while convergence to critical points (points with vanishing gradients) is *not*. Theoretical works have also affirmed the same phenomenon for other machine learning problems (see [5, 6, 18, 19] and the references therein).

In this paper we give a provable linear-time algorithm for finding an *approximate* local minimum in smooth non-convex optimization. It applies to a general setting of machine learning optimization, and in particular to the optimization problem of training deep neural networks. Furthermore, the running time bound of our algorithm is the fastest known even for the more lenient task of computing a point with vanishing gradient (called a critical point), for a wide range of parameters.

Formally, the problem of unconstrained mathematical optimization is stated in general terms as that of finding the minimum value that a function attains over Euclidean space, i.e.

$$\min_{x \in \mathbb{R}^d} f(x) \ . \tag{1.1}$$

If $f$ is convex, the above formulation is *convex optimization* and is solvable in (randomized) polynomial time even if only a valuation oracle to $f$ is provided. A crucial property of convex functions is that "local optimality implies global optimality", allowing for greedy algorithms to reach the global optimum efficiently. Unfortunately, this is no longer the case if $f$ is nonconvex; indeed, even a degree four polynomial can be NP-hard to optimize [23], or even just to check whether a point is not a local minimum [25]. Thus, for non-convex optimization one has to settle for the more modest goal of reaching approximate local optimality efficiently.

Note that a particular interest to machine learning is the optimization of functions $f : \mathbb{R}^d \mapsto \mathbb{R}$ of the finite-sum form

$$f(x) = \frac{1}{n} \sum_{i=1}^{n} f_i(x) \ . \tag{1.2}$$

Such functions arise when minimizing loss over a training set, where each example $i$ in the set corresponds to one loss function $f_i$ in the summation.

We say that the function $f$ is second-order smooth if it has Lipschitz continuous gradient and Lipschitz continuous Hessian. We say that a point $x$ is an $\varepsilon$-approximate local minimum if it satisfies (following the tradition of [28]):

$$\|\nabla f(x)\| \leq \varepsilon \quad \text{and} \quad \nabla^2 f(x) \succeq -\sqrt{\varepsilon}\mathbf{I} \ ,$$

where $\|\cdot\|$ denotes the Euclidean norm of a vector. We say that a point $x$ is an $\varepsilon$-critical point if it satisfies the gradient condition above, but not necessarily the second-order condition. Critical points include saddle points in addition to local minima. We remark that $\varepsilon$-approximate local minima (even with $\varepsilon = 0$) are not necessarily close to any local minimum, neither in domain nor in function value. However, if we assume in addition the function satisfies the (robust) strict-saddle property [15, 24] (see Section 2 for the precise definition), then an $\varepsilon$-approximate local minimum is guaranteed to be close to a local minimum for sufficiently small $\varepsilon$.

Our main theorem below states the time required for the proposed algorithm FastCubic to find an $\varepsilon$-approximate local minimum for second-order smooth functions.

**Theorem** (informal). *Ignoring smoothness parameters, the running time of* FastCubic *to return an $\varepsilon$-approximate local minimum is*

$$\tilde{O}\left(\left(\frac{n}{\varepsilon^{3/2}} + \frac{n^{3/4}}{\varepsilon^{7/4}}\right) \cdot \mathbb{T}_{h,1}\right) \text{ for (1.2)} \quad \text{or} \quad \tilde{O}\left(\frac{1}{\varepsilon^{7/4}} \cdot \mathbb{T}_h\right) \text{ for the general (1.1).}$$

*Above, $\mathbb{T}_h$ is the time to compute Hessian-vector product for $\nabla^2 f(x)$ and $\mathbb{T}_{h,1}$ is that for an arbitrary $\nabla^2 f_i(x)$.*

The full statement of Theorem 1 can be found in Section 2.

Hessian-vector products can be computed in linear time —meaning $\mathbb{T}_{h,1} = O(d)$ and $\mathbb{T}_h = O(nd)$— for many machine learning problems such as generalized linear models and training neural networks [1, 29]. We explain this more generally in Appendix A. Therefore,

**Corollary 1.1.** *Algorithm* FastCubic *returns an $\varepsilon$-approximate local minimum for the optimization problem of training a neural network in time*

$$\tilde{O}\left(\frac{nd}{\varepsilon^{3/2}} + \frac{n^{3/4}d}{\varepsilon^{7/4}}\right) \ .$$

Another important aspect of our algorithm is that even in terms of just reaching an $\varepsilon$-critical point, i.e. a point that satisfies $\|\nabla f(x)\| \leq \varepsilon$ without any second-order guarantee, FastCubic is faster than all previous results (see Table 1 for a comparison).

The fastest methods to find critical points for a smooth non-convex function are gradient descent and its extensions, jointly known as first-order methods. These methods are extremely efficient in terms of per-iteration complexity; however, they necessarily suffer from a $1/\varepsilon^2$ convergence rate [27], to the best of our knowledge, in previous results only higher-order methods seem capable of breaking this $1/\varepsilon^2$ bottleneck [28]. For certain ranges of parameters, our FastCubic finds local minima even faster than first-order methods, even though they only find critical points. This is depicted in Table 1.

| Paper | Total Time Achieving $\|\nabla f(x)\| \leq \varepsilon$ | Second-Order Guarantee |
|---|---|---|
| Gradient Descent (GD) | $O\left(\frac{nd}{\varepsilon^2}\right)$ | n/a |
| SVRG [2] | $O\left(nd + \frac{n^{2/3}d}{\varepsilon^2}\right)$ | n/a |
| SGD [20] | $O\left(\frac{d}{\varepsilon^4}\right)$ | n/a |
| noisy SGD [16] [a] | $O\left(\frac{d^{C_1}}{\varepsilon^4}\right)$ | $\nabla^2 f(x) \succeq -\varepsilon^{1/C_2}\mathbf{I}$ |
| cubic regularization [28] | $\tilde{O}\left(\frac{nd^{\omega-1}+d^\omega}{\varepsilon^{3/2}}\right)$ | $\nabla^2 f(x) \succeq -\varepsilon^{1/2}\mathbf{I}$ |
| **this paper** | $\tilde{O}\left(\frac{nd}{\varepsilon^{7/4}}\right)$ | $\nabla^2 f(x) \succeq -\varepsilon^{1/2}\mathbf{I}$ |
| **this paper** | $\tilde{O}\left(\frac{nd}{\varepsilon^{3/2}} + \frac{n^{3/4}d}{\varepsilon^{7/4}}\right) $ | $\nabla^2 f(x) \succeq -\varepsilon^{1/2}\mathbf{I}$ |

Table 1: Comparison of known methods.

[a] Here $C_1, C_2$ are two constants that are not explicitly written. We believe $C_1 \geq 4$.

## 1.1 Related work

**Methods that Provably Reach Critical Points.** Recall that only a gradient oracle is needed to reach a critical point. The most commonly used algorithm in practice for training non-convex learning machines such as deep neural networks is stochastic gradient descent (SGD), also known as stochastic approximation [30] and its derivatives. Some practical enhancements widely used in practice are based on Nesterov's acceleration [26] and adaptive regularization [12]. The variance reduction technique, introduced in [32], was extremely successful in convex optimization, but only recently there was a non-convex counterpart with theoretical benefits introduced [2].

**Methods that Provably Reach Local Minima.** The recent work of Ge *et al.* [17] showed that a noise-injected version of SGD in fact converges to local minima instead of critical points, as long as the underlying non-convex function is strict-saddle. Their theoretical running time is a large polynomial in the dimension and not competitive with our method (see Table 1).

The work of Lee et al. [24] shows that gradient descent, starting from a random point, almost surely converges to a local minimum of a strict-saddle function. The rates of convergence and precise step-sizes that are required are, however, yet unknown.

If second-order information (i.e., the Hessian oracle) is provided, the cubic-regularization method of Nesterov and Polyak [28] converges in $O(\frac{1}{\varepsilon^{3/2}})$ iterations. However, each iteration of Nesterov-Polyak requires solving a cubic function which, in general, takes time super-linear in the input representation.

One natural direction is to apply an approximate trust region solver, such as the linear-time solver of [22], to approximately solve the cubic regularization subroutine of Nesterov-Polyak. However, the approximation needed by a naive calculation makes this approach even slower than vanilla gradient descent. Our main challenge is to obtain approximate second-order local-minima and simultaneously improve upon gradient descent.

Independently of this paper and concurrently[1], Carmon *et al.* [7] develop an accelerated gradient descent method that achieves the same running time for finding an approximate local minimum as in our paper. Remarkably, the same running time is obtained via a very different technique.

## 1.2 Our Techniques

Our algorithm is based on the cubic regularization method of Nesterov and Polyak [8, 9, 28]. At a high level, cubic regularization states that if we can minimize a cubic function $m(h) \triangleq g^\top h + \frac{1}{2}h^\top \mathbf{H}h + \frac{L}{6}\|h\|^3$ exactly, where $g = \nabla f(x)$, $\mathbf{H} = \nabla^2 f(x)$, and $L$ is the second-order smoothness of the function $f$, then we can iteratively perform updates $x' \leftarrow x + h$, and this algorithm converges to an $\varepsilon$-approximate local minimum in $O(1/\varepsilon^{3/2})$ iterations. Unfortunately, solving this cubic minimization problem exactly, to the best of our knowledge, requires a running time of $O(d^\omega)$ where $\omega$ is the matrix multiplication constant. Getting around this requires five observations.

The *first observation* is that, minimizing $m(h)$ up to a constant multiplicative approximation (plus a few other constraints) is sufficient for showing an iteration complexity of $O(1/\varepsilon^{3/2})$.[2] The proof techniques to show this observation are based on extending Nesterov and Polyak.

The *second observation* is that the minimizer $h^*$ of $m(h)$ must be of the form $h^* = (\mathbf{H}+\lambda^*\mathbf{I})^+g + v$, where $\lambda^* \geq 0$ is some constant satisfying $\mathbf{H} + \lambda^*\mathbf{I} \succeq 0$, and $v$ is the smallest eigenvector of $\mathbf{H}$ and $^+$ denotes the pseudo-inverse of a matrix. This can be viewed as moving in a mixture direction between choosing $h \leftarrow v$, and choosing $h$ to follow a shifted Newton's direction $h \leftarrow (\mathbf{H} + \lambda^* \mathbf{I})^+ g$. Intuitively, we wish to reduce both the computation of $(\mathbf{H} + \lambda^* \mathbf{I})^+ g$ and $v$ to Hessian-vector products.

---

[1] To be precise, their manuscript appeared online approximately 24 hours before ours.
[2] More specifically, we need $m_t(h) \leq \frac{1}{C}\min_h\{m_t(h)\}$ for some constant $C$. In addition, we need to have good bounds on $\|h\|$ and $\|\nabla m(h)\|$.

The first task of computing $(\mathbf{H} + \lambda^* \mathbf{I})^+ g$ can be slow, and even if $\mathbf{H} + \lambda^* \mathbf{I}$ is strictly positive-definite, computing it has a complexity depending on the (possibly huge) condition number of $\mathbf{H} + \lambda^* \mathbf{I}$ [34]. The *third observation* is that it suffices to pick some $\lambda' > \lambda^*$ so both (1) the condition number of $\mathbf{H} + \lambda' \mathbf{I}$ is small and (2) the vectors $(\mathbf{H} + \lambda^* \mathbf{I})^{-1} g$ and $(\mathbf{H} + \lambda' \mathbf{I})^{-1} g$ are close. This relies on the structure of $m(h)$.

The second task of computing $v$ has a complexity depending on $1/\sqrt{\delta}$ where $\delta$ is the target additive error [13, 14]. The *fourth observation* is that the choice $\delta = \sqrt{\varepsilon}$ suffices for the outer loop of cubic regularization to make sufficient progress. This reduces the complexity to compute $v$.

Finally, finding the correct value $\lambda^*$ itself is as hard as minimizing $m_t(h)$. The *fifth step* is to design an iterative scheme that makes only logarithmic number of guesses on $\lambda^*$. This procedure either finds the correct one (via binary search), or finds an approximate one, $\lambda'$, but satisfying $(\mathbf{H} + \lambda^* \mathbf{I})^{-1} g$ and $(\mathbf{H} + \lambda' \mathbf{I})^{-1} g$ being sufficiently close.

Putting all the observations together, and balancing all the parameters, we can obtain a cubic minimization subroutine (see FastCubicMin in Algorithm 2) that runs in time $O(nd + n^{3/4} d/\varepsilon^{1/4})$.

## 2 Preliminaries and Main Theorem

We use $\|\cdot\|$ to denote the Euclidean norm of a vector and the spectral norm of a matrix. For a symmetric matrix $\mathbf{M}$ we denote by $\lambda_{\max}(\mathbf{M})$ and $\lambda_{\min}(\mathbf{M})$ respectively the maximum and minimum eigenvalues of $\mathbf{M}$. We denote by $\mathbf{A} \succeq \mathbf{B}$ that $\mathbf{A} - \mathbf{B}$ is positive semidefinite (PSD). For a PSD matrix $\mathbf{M}$, we denote by $\mathbf{M}^+$ its pseudo-inverse if $\mathbf{M}$ is not strictly positive definite.

We make the following Lipschitz continuity assumptions for the gradient and Hessian of the target function $f$. Namely, there exist $L_2, L > 0$ such that

$$\forall x, y \in \mathbb{R}^d: \quad \|\nabla^2 f(x)\| \leq L_2 \quad \text{and} \quad \|\nabla^2 f(x) - \nabla^2 f(y)\| \leq L\|x - y\|. \tag{2.1}$$

**Definition 2.1.** *We assume the following complexity parameters on the access to $f(x)$:*

- *Let $\mathbb{T}_g \in \mathbb{R}^*$ be the time complexity to compute $\nabla f(x)$ for any $x \in \mathbb{R}^d$.*
- *Let $\mathbb{T}_h \in \mathbb{R}^*$ be the time complexity to compute $(\nabla^2 f(x)) v$ for any $x, v \in \mathbb{R}^d$.*

**Definition 2.2.** *We say that $f$ is of finite-sum form if $f = \frac{1}{n} \sum_{i=1}^n f_i(x)$ and $\|\nabla^2 f_i(x)\| \leq L_2$ for each $i \in [n]$. In this case, we define $\mathbb{T}_{h,1}$ to be the time complexity to compute $(\nabla^2 f_i(x)) v$ for arbitrary $x, v \in \mathbb{R}^d$ and $i \in [n]$.*

Next we define the strict-saddle function for which an $\varepsilon$-approximate local minimum is almost equivalent to a local minimum [15, 24].

**Definition 2.3** (strict saddle). *Suppose $f(\cdot): \mathbb{R}^d \to \mathbb{R}$ is twice differentiable. For $\alpha, \beta, \gamma \geq 0$, we say $f$ is $(\alpha, \beta, \gamma)$-strict saddle if every $x \in \mathbb{R}^d$ satisfies at least one of the following three conditions:*
1. $\|\nabla f(x)\| \geq \alpha$.
2. $\lambda_{\min}(\nabla^2 f) \leq -\beta$.
3. *There exists a local minimum $x^\star$ that is $\gamma$-close to $x$ in Euclidean distance.*

We see that if a function is $(\alpha, \beta, \gamma)$-strict saddle, then for $\varepsilon < \min\{\alpha, \beta^2\}$ an $\varepsilon$-approximate local minimum is $\gamma$-close to some local minimum.

---
**Algorithm 1** FastCubic($f, x_0, \varepsilon, L, L_2$)
---
**Input:** $f(x)$ that satisfies (2.1) with $L_2$ and $L$; a starting vector $x_0$; a target accuracy $\varepsilon$.
1: $\kappa \leftarrow \left(\frac{900}{\varepsilon L}\right)^{1/2}$.
2: **for** $t = 0$ **to** $\infty$ **do**
3: $\quad m_t(h) \triangleq \nabla f(x_t)^\top h + \frac{h^\top \nabla^2 f(x_t) h}{2} + \frac{L}{6}\|h\|^3$
4: $\quad (\lambda, v, v_{\min}) \leftarrow \mathsf{FastCubicMin}\left(\nabla f(x_t), \nabla^2 f(x_t), L, L_2, \kappa\right)$
5: $\quad h' \leftarrow$ either $v$ or $\frac{\lambda v_{\min}}{2L}$ whichever gives smaller value for $m_t(h)$;
6: $\quad$ Set $x_{t+1} \triangleq x_t + h'$
7: $\quad$ **if** $m_t(h') > -\frac{\varepsilon^{3/2}}{c\sqrt{L}}$ **then return** $x_{t+1}$.  $\diamond$ *c is a constant; we proved $c = 2.4 * 10^6$ works*
8: **end for**
---

## 2.1 Main Results

The finite-sum setting captures much of supervised learning, including Neural Networks and Generalized Linear Models. The main theorem which we show in our paper is as follows:

> **Theorem 1.** FastCubic *(Algorithm 1) starts from a point $x_0$ and outputs a point $x$ such that*
> $$\|\nabla f(x)\| \leq \varepsilon \quad \text{and} \quad \lambda_{\min}(\nabla^2 f(x)) \geq -\sqrt{L\varepsilon}$$
> *in total time (denoting by $D \triangleq f(x_0) - f(x^*)$)*
> - $\tilde{O}\left(\frac{D\sqrt{L}}{\varepsilon^{3/2}} \cdot \mathbb{T}_g + \frac{DL^{1/4}\sqrt{L_2}}{\varepsilon^{7/4}} \cdot \mathbb{T}_h\right)$, *or*
> - $\tilde{O}\left(\frac{D\sqrt{L}}{\varepsilon^{3/2}} \cdot (\mathbb{T}_g + n\mathbb{T}_{h,1}) + \frac{Dn^{3/4}L^{1/4}\sqrt{L_2}}{\varepsilon^{7/4}} \cdot \mathbb{T}_{h,1}\right)$ *in the finite-sum setting (see Definition 2.2).*
>
> *Here $\tilde{O}$ hides logarithmic factors in $L, L_2, 1/\varepsilon, d$, and in $\max_x \{\|\nabla f(x)\|\}$.*

**Two Known Subroutines.** Our running time of FastCubic relies on the following recent results for approximate matrix inverse and approximate PCA:

**Theorem 2.4** (Approximate Matrix Inverse). *Suppose matrix $\mathbf{M} \in \mathbb{R}^{d \times d}$ satisfies $\|\mathbf{M}\| \leq L_2$ and $\lambda \mathbf{I} + \mathbf{M} \succeq \delta \mathbf{I}$ for constants $\lambda, \delta, L_2 > 0$. Let $\kappa \triangleq \frac{\lambda + L_2}{\delta}$. Then, we can compute vector $x$ satisfying*
$$\|x - (\lambda \mathbf{I} + \mathbf{M})^{-1} b\| \leq \varepsilon \|b\|, \tag{2.2}$$
*using Accelerated gradient descent (AGD) in $O(\kappa^{1/2} \log(\kappa/\varepsilon))$ iterations, each requiring $O(d)$ time plus the time needed to multiply $\mathbf{M}$ with a vector.*

*Moreover, suppose $\mathbf{M} = \frac{1}{n}\sum_{i=1}^n \mathbf{M}_i$ where each $\mathbf{M}_i$ is symmetric and satisfies $\|\mathbf{M}_i\| \leq L_2$. If $\mathbf{M}_i b$ can be computed in time $O(d')$ for each $i$ and vector $b$, then accelerated SVRG [4, 33] computes a vector $x$ that satisfies equation (2.2) in time $O\left(\max\{n, n^{3/4}\kappa^{1/2}\} \cdot d' \cdot \log^2(\kappa/\varepsilon)\right)$.*

*We refer to the running time for this computation as $\mathbb{T}_{\mathsf{inverse}}(\kappa, \varepsilon)$ and the algorithm as $\mathcal{A}$.*

Above, the SVRG based running time shall be used only towards our finite-sum case in Definition 2.2.

**Theorem 2.5** (AppxPCA [3, 13, 14]). *Let $\mathbf{M} \in \mathbb{R}^{d \times d}$ be a symmetric matrix with eigenvalues $1 \geq \lambda_1 \geq \cdots \geq \lambda_d \geq 0$. With probability at least $1-p$, AppxPCA produces a unit vector $w$ satisfying $w^\top M w \geq (1-\delta_\times)(1-\varepsilon)\lambda_{\max}(\mathbf{M})$. The total running time is $\tilde{O}(\mathbb{T}_{\mathsf{inverse}}(1/\delta_\times, \varepsilon\delta_\times))$.*

## 3 Our Fast Cubic Regularization Algorithm

Recall that the cubic regularization method of Nesterov and Polyak [28] studies the following upper bound on the change in objective value as we move from a point $x_t$ to $x_t + h$: (it follows simply from the Taylor series truncated to the third order)

$$\forall h \in \mathbb{R}^d: \quad f(x_t + h) - f(x_t) \leq m_t(h) \triangleq \nabla f(x_t)^\top h + \frac{h^\top \nabla^2 f(x_t) h}{2} + \frac{L}{6} \|h\|^3 \ . \qquad (3.1)$$

Denote by $h^*$ an arbitrary minimizer of $m_t(h)$. We propose in this paper a subroutine FastCubicMin to minimizes $m_t(h)$ approximately. Note that FastCubicMin returns two vectors $v$ and $v_{\min}$. We then choose $h'$ to be either $v$ or $\frac{\lambda v_{\min}}{2L}$, whichever gives a smaller value for $m_t(h)$.

Before discussing the details of FastCubicMin, let us first state a main theorem for FastCubicMin:[3]

---

**Theorem 2** (Guarantees of FastCubicMin). *The algorithm FastCubicMin finds a vector $h'$ that satisfies*

(a) *It produces a vector $h'$ satisfying $m_t(h') \leq 0$ and*

$$\text{either} \quad 3000 m_t(h') \leq m_t(h^*) \quad \text{or} \quad m_t(h^*) \geq -\frac{\varepsilon^{3/2}}{800\sqrt{L}} \ .$$

(b) *If $m(h^*) \geq -\frac{\varepsilon^{3/2}}{300\sqrt{L}}$, then $\|h'\| \leq \|h^*\| + \frac{\sqrt{\varepsilon}}{4\sqrt{L}}$ and $\|\nabla m_t(h')\| \leq \frac{\varepsilon}{2}$.*

(c) FastCubicMin *runs in time: (using $\tilde{O}$ to hide logarithmic factors in $L, L_2, 1/\varepsilon, d, \|\nabla f(x_t)\|$)*

- $\tilde{O}\left(\frac{\sqrt{L_2}}{(\varepsilon L)^{1/4}} \cdot \mathbb{T}_h\right)$ *where $\mathbb{T}_h$ is the time to multiply $\nabla^2 f(x_t)$ to a vector;*
- $\tilde{O}\left(\max\left\{n, n^{3/4} \frac{\sqrt{L_2}}{(\varepsilon L)^{1/4}}\right\} \cdot \mathbb{T}_{h,1}\right)$ *where $\mathbb{T}_{h,1}$ is the time to multiply $\nabla^2 f_i(x_t)$ with a vector.*

---

Above, the first guarantee promises that we are either done (because $m_t(h^*)$ is close to zero), or we obtain a 1/3000 *multiplicative* approximation to $m_t(h^*)$. Our second guarantee in Theorem 2 promises that when we are done (because $m_t(h^*)$ is close to zero), the output vector $h'$ and $h^*$ are roughly similar in Euclidean norm and have a small gradient $\|\nabla m_t(h')\|$. Our third guarantee gives the time complexity of FastCubicMin.

Now, our final algorithm FastCubic for finding the $\varepsilon$-approximate local minimum of $f(x)$ is included in Algorithm 1. It simply iteratively calls FastCubicMin to find an approximate minimizer, and it then stops whenever $m_t(h') > -\frac{\varepsilon^{3/2}}{c\sqrt{L}}$ for some large constant $c$.

**Roadmap.** In Section 4 we show why Theorem 2 implies Theorem 1. All the remaining sections are for the purpose of proving Theorem 2. Because our FastCubicMin is very technical, instead of stating what the algorithm is right away, we decide to take a different path. In Section 5, we first state a lemma characterizing "what $h^*$ looks like". In Section 6, we provide a set of sufficient conditions which "look similar" to the characterization of $h^*$, and show that as long as these conditions are met, Theorem 2-a and 2-b follow easily. Finally, in Section 7, we state FastCubicMin and explain why it satisfies these sufficient conditions and why it runs in the aforementioned time.

## 4 Theorem 2 implies Theorem 1

In this section, we show that Theorem 2 implies Theorem 1. It relies on the following lemma (proved in Appendix B) regarding the sufficient condition for us to reach an $\varepsilon$-approximate local minimum.

---

[3]To present the simplest result, we have not tried to improve the constant dependency in this paper.

**Lemma 4.1.** *If $m_t(h^*) \geq -\frac{\varepsilon^{3/2}}{800\sqrt{L}}$ and $h'$ is an approximate minimizer of $m_t(h)$ satisfying*

$$\|h'\| \leq \|h^*\| + \frac{\sqrt{\varepsilon}}{4\sqrt{L}} \quad \text{and} \quad \|\nabla m_t(h')\| \leq \frac{\varepsilon}{2},$$

*then we have that $\|\nabla f(x_t + h')\| \leq \varepsilon$ and $\lambda_{\min}(\nabla^2 f(x_t + h')) \geq -\sqrt{L\varepsilon}$.*

*Proof of Theorem 1 from Theorem 2.* When FastCubic terminates, we have $m_t(h') > -\frac{\varepsilon^{3/2}}{c\sqrt{L}}$; therefore, it satisfies $m_t(h^*) \geq -\frac{\varepsilon^{3/2}}{800\sqrt{L}}$ according to Theorem 2-a. Combining this with Theorem 2-b and Corollary 4.1, we conclude that in the last iteration of FastCubic, our output satisfies $\|\nabla f(x_t+h')\| \leq \varepsilon$ and $\lambda_{\min}(\nabla^2 f(x_t+h')) \geq -\sqrt{L\varepsilon}$. This finishes the proof with respect to the accuracy conditions.

As for the running time, in every iteration except for the last one, FastCubic satisfies $m_t(h') \leq -\Omega\left(\frac{\varepsilon^{3/2}}{\sqrt{L}}\right)$. Therefore by (3.1), we must have decreased the objective by at least $\Omega\left(\frac{-\varepsilon^{3/2}}{\sqrt{L}}\right)$ in this round, and this cannot happen for more than $O\left(\frac{(f(x_0) - f^*)\sqrt{L}}{\varepsilon^{3/2}}\right)$ iterations. The final running time of FastCubic follows from this bound together with Theorem 2-c. □

Therefore, in the rest of the paper it suffices to study FastCubicMin and prove Theorem 2.

## 5 Characterization Lemma of the Minimizer $h^*$

For notational simplicity in this and the subsequent sections we focus on the following problem:

$$\text{minimize} \quad m(h) \triangleq g^\top h + \frac{h^\top \mathbf{H} h}{2} + \frac{L}{6}\|h\|^3$$

$$\text{where } \mathbf{H} \text{ is a symmetric matrix with } \|\mathbf{H}\|_2 \leq L_2.$$

Recall from the previous section that we have denoted by $h^*$ an arbitrary minimizer of $m(h)$. We have the following lemma which characterizes $h^*$: (a variant of this lemma has appeared in [8], and we prove it in the appendix for the sake of completeness)

**Lemma 5.1.** *We have $h^*$ is a minimizer of $m(h)$ if and only if there exists $\lambda^* \geq 0$ such that*

$$\mathbf{H} + \lambda^*\mathbf{I} \succeq 0, \quad (\mathbf{H} + \lambda^*\mathbf{I})h^* = -g, \quad \|h^*\| = \frac{2\lambda^*}{L}.$$

*The objective value in this case is given by*

$$m(h^*) = -\frac{1}{2}g^\top(\mathbf{H} + \lambda^*\mathbf{I})^+ g - \frac{2(\lambda^*)^3}{3L^2} \leq 0.$$

The following corollary comes from Lemma 5.1 and its proof:

**Corollary 5.2.** *The value $\lambda^*$ in Lemma 5.1 is unique, and for every $\lambda$ satisfying $\mathbf{H} + \lambda\mathbf{I} \succ 0$, we have*

$$\|(\mathbf{H} + \lambda\mathbf{I})^{-1}g\| > \frac{2\lambda}{L} \iff \lambda^* > \lambda \quad \text{and} \quad \|(\mathbf{H} + \lambda\mathbf{I})^{-1}g\| < \frac{2\lambda}{L} \iff \lambda^* < \lambda.$$

In the above characterization, we have a crude upper bound on $\lambda^*$:

**Proposition 5.3.** *We have $\lambda^* \leq B \triangleq \max\{2L_2 + \sqrt{L\|g\|}, 1\}$ with $\lambda^*$ defined in Lemma 5.1.*

*Proof.* We have $L\|(\mathbf{H} + B\mathbf{I})^{-1}g\| \leq \frac{L\|g\|}{\lambda_{\min}(\mathbf{H}+B\mathbf{I})} \leq \frac{L\|g\|}{B-L_2} < 2B$ and therefore $\lambda^* \leq B$ due to Corollary 5.2. □

# 6 Sufficient Conditions for Theorem 2-a and 2-b

Without worrying about the design of FastCubicMin at this moment, let us first state a set of sufficient conditions under which the assumptions in Theorem 2-a can be satisfied.

**Main Lemma 1.** *Consider an algorithm that outputs a real $\lambda \in [0, 2B]$, a vector $v \in \mathbb{R}^d$, and a unit vector $v_{\min} \in \mathbb{R}^d$. Additionally, suppose numbers $\kappa, \tilde{\varepsilon} \geq 0$ satisfying the following conditions:*

$$\tilde{\varepsilon} \leq \frac{1}{10000} \frac{1}{(\max\{\kappa, L, L_2, \|g\|, \|(\mathbf{H} + \lambda \mathbf{I})^{-1}\|, B\})^{20}} \tag{6.1}$$

$$(\mathbf{H} + (\lambda - L\tilde{\varepsilon})\mathbf{I})^{-1} \succ 0 \tag{6.2}$$

*Moreover, suppose that the outputs $(\lambda, v, v_{\min})$ satisfy one of the following two cases:*

**Case 1:** $L\|(\mathbf{H} + \lambda \mathbf{I})^{-1} g\| \in [2\lambda - 2L\tilde{\varepsilon}, 2\lambda + 2L\tilde{\varepsilon}]$ *and* $\|v + (\mathbf{H} + \lambda \mathbf{I})^{-1} g\| \leq \tilde{\varepsilon}$

**Case 2:** *The following conditions are satisfied:*

(a) $\lambda \geq \lambda^*$ *and* $\lambda + \lambda_{\min}(\mathbf{H}) \leq \frac{1}{\kappa}$

(b) $L\|(\mathbf{H} + \lambda \mathbf{I})^{-1} g\| \leq 2\lambda$ *and* $\|v + (\mathbf{H} + \lambda \mathbf{I})^{-1} g\| \leq \tilde{\varepsilon}$

(c) $v_{\min}^\top \mathbf{H} v_{\min} \leq \lambda_{\min}(\mathbf{H}) + \frac{1}{10\kappa}$

*Then, at least one of the two choices $h' \in \{v, \frac{\lambda v_{\min}}{2L}\}$ satisfy*

$$\text{either} \quad m(h^*) \geq 3000 m(h') \quad \text{or} \quad m(h^*) \geq -\frac{32}{\kappa^3 L^2} \ .$$

Let us compare such sufficient conditions to the characterization Lemma 5.1.

- In Case 1, up to a very small error $\tilde{\varepsilon}$, we have essentially found a vector $v$ that satisfies $v \approx -(\mathbf{H} + \lambda \mathbf{I})^{-1} g$ and $\|v\| \approx \frac{2\lambda}{L}$. Therefore, this $v$ should be close to $h^*$ for obvious reason. (This is the simple case.)

- In Case 2, we have only found a vector $v$ that satisfies $v \approx -(\mathbf{H} + \lambda \mathbf{I})^{-1} g$ and $\|v\| \lesssim \frac{2\lambda}{L}$. In this case, we also compute an approximate lowest eigenvector $v_{\min}$ of $\lambda_{\min}(\mathbf{H})$ up to an additive $1/10\kappa$ accuracy (see case 2-c). We will make sure that, as long as the conditions in 2-a hold, then either $v$ or $\frac{\lambda v_{\min}}{2L}$ will be an approximate minimizer for $m_t(h)$.

    (This is the hard case.)

*Proof of Main Lemma 1.* We first consider Case 1. According to Corollary 5.2, if $\tilde{\varepsilon} = 0$ then $v$ is a minimizer of $m(h)$. The following claim extends this argument to the setting when $\tilde{\varepsilon} > 0$:

**Claim 6.1.** *If $\lambda$ and $v$ satisfy Case 1 and $\tilde{\varepsilon}$ satisfies (6.1), then $m(v) \leq m(h^*) + \frac{1}{250\kappa^3 L^2}$*

From the above lemma it follows that either $m(h^*) \geq -\frac{8}{\kappa^3 L^2}$ otherwise $m(h^*) \geq 1.1 m(v)$ which satisfies the conditions of the theorem.

We now consider Case 2, and in this case we make the following two claims:

**Claim 6.2.** *If $\lambda_{\min}(\mathbf{H}) \leq -\frac{1}{\kappa}$ then $m(h^*) \geq 1500 \min\left\{m(v), m\left(\frac{\lambda v_{\min}}{2L}\right)\right\} - \frac{1}{500\kappa^3 L^2}$ .*

**Claim 6.3.** *If $\lambda_{\min}(\mathbf{H}) \geq -\frac{1}{\kappa}$ then $m(h^*) \geq 2m(v) - \frac{16}{\kappa^3 L^2}$ .*

Lemma 1 now follows from the two claims because we can output the vector $h'$ which has the lowest value of $m(h')$ amongst the two choices $h' \in \{v, \lambda \frac{v_{\min}}{2d}\}$. This satisfies either $m(h^*) \geq 3000 m(h')$ or $m(h^*) \geq -\frac{32}{\kappa^3 L^2}$.

The missing proofs of the three claims are deferred to Appendix D. □

The next main lemma shows that, under the same sufficient conditions as Main Lemma 1, we also have that Theorem 2-b holds. (Its proof is contained in Appendix E.)

**Main Lemma 2.** *In the same setting as Main Lemma 1, suppose $m(h^*) \geq -\frac{\varepsilon^{3/2}}{300\sqrt{L}}$. Then the output vector $v$ satisfies the following conditions:*
$$\|v\| \leq \|h^*\| + \frac{3}{\kappa L} \quad \text{and} \quad \|\nabla m(v)\| \leq \frac{\varepsilon}{4} + \frac{15}{\kappa^2 L} .$$

# 7 Main Algorithms for Theorem 2

We are now ready to state our main algorithm FastCubicMin and sketch why it satisfies the sufficient conditions in Main Lemma 1. As described in Algorithm 2, our algorithm starts with a very large choice $\lambda_0 \leftarrow 2B$ and decreases it gradually. At each iteration $i$, it computes an approximate inverse $v$ satisfying $\|v + (\mathbf{H} + \lambda_i \mathbf{I})^{-1} g\| \leq \tilde{\varepsilon}$ with respect to the current $\lambda_i$. Then there are three cases, depending on whether $L\|v\|$ is approximately equal to, larger than, or smaller than $2\lambda_i$. At a high level, if it is "equal", then we have met Case 1 in Main Lemma 1; if it is "larger", then we can binary search the correct value of $\lambda^*$ in the interval $[\lambda_i, \lambda_{i-1}]$; and if it is "smaller", then we need to compute an approximate eigenvector and carefully choose the next point $\lambda_{i+1}$.

We state our main lemma below regarding the correctness and running time of FastCubicMin.

**Main Lemma 3.** FastCubicMin *in Algorithm 2 outputs a real $\lambda \in [0, 2B]$, a vector $v \in \mathbb{R}^d$, and a unit vector $v_{\min} \in \mathbb{R}^d$ satisfying one of the two sufficient conditions in Main Lemma 1. We also have that the procedure can be implemented in a total running time of*

- $\tilde{O}(\sqrt{\kappa L_2} \cdot \mathbb{T}_h)$ *if Accelerated Gradient Descent is used in Theorem 2.4 to invert matrices.*
- $\tilde{O}(\max\{n, n^{3/4}\sqrt{\kappa L_2}\} \cdot \mathbb{T}_{h,1})$ *if we use accelerated SVRG as the subprocedure $\mathcal{A}$ in Theorem 2.4.*

*Here $\tilde{O}$ hides logarithmic factors in $L, L_2, \kappa, d, B$.*

We prove the correctness half of Main Lemma 3, and defer its running time analysis to Appendix G.

## 7.1 Correctness Half of Main Lemma 3

We will now establish the correctness of our algorithm. We first observe that the BinarySearch subroutine returns $(\lambda, v, \emptyset)$ that satisfies Case 1 of Main Lemma 1.

**Fact 7.1.** BinarySearch *outputs a pair $\lambda$ and $v$ such that*
$$L\|(\mathbf{H} + \lambda\mathbf{I})^{-1}g\| \in [2\lambda - 2L\tilde{\varepsilon}, 2\lambda + 2L\tilde{\varepsilon}] \quad \text{and} \quad \|v + (\mathbf{H} + \lambda\mathbf{I})^{-1}g\| \leq \tilde{\varepsilon} .$$

*Proof.* The latter is guaranteed by line 3 in BinarySearch, and the former is implied by the latter because
$$L\|(\mathbf{H} + \lambda\mathbf{I})^{-1}g\| \in \left[L\|v\| - L\tilde{\varepsilon}/2, L\|v\| + L\tilde{\varepsilon}/2\right] \subseteq \left[2\lambda - 2L\tilde{\varepsilon}, 2\lambda + 2L\tilde{\varepsilon}\right] . \qquad \square$$

We also establish the following invariants regarding the values $\lambda_i$. (Proof in Appendix F.)

**Lemma 7.2.** *The following statements hold for all $i$ until* FastCubicMin *terminates*

(a) $\lambda_i \in [0, 2B]$, $\lambda_i + \lambda_{\max}(\mathbf{H}) \leq 3B$

(b) $\lambda_i + \lambda_{\min}(\mathbf{H}) \geq \frac{3}{10\kappa}$

(c) $\lambda_{i+1} + \lambda_{\min}(\mathbf{H}) \leq \frac{3}{4}(\lambda_i + \lambda_{\min}(\mathbf{H}))$ *unless $\lambda_{i+1} = 0$*

*Moreover when* FastCubicMin *terminates at Line 20 we have $\lambda_i + \lambda_{\min}(\mathbf{H}) \leq \frac{1}{\kappa}$.*

We now prove the output $(\lambda, v, v_{\min})$ of FastCubicMin satisfies the sufficient conditions of Main Lemma 1.

**Algorithm 2** FastCubicMin($g, \mathbf{H}, L, L_2, \kappa$) (main algorithm for cubic minimization)

**Input:** $g$ a vector, $\mathbf{H}$ a symmetric matrix, parameters $\kappa, L$ and $L_2$ which satisfies $-L_2\mathbf{I} \preceq \mathbf{H} \preceq L_2\mathbf{I}$.
**Output:** $(\lambda, v, v_{\min})$
1: $B \leftarrow L_2 + \sqrt{L\|g\|} + \frac{1}{\kappa}$.
2: $\tilde{\varepsilon} \leftarrow 1/\left(10000 \left(\max\left\{L, \|g\|, \frac{3\kappa}{10}, B, 1\right\}\right)^{20}\right)$
3: $\lambda_0 \leftarrow 2B$.
4: **for** $i = 0$ **to** $\infty$ **do**
5:     Compute $v$ such that $\|v + (\mathbf{H} + \lambda_i\mathbf{I})^{-1}g\| \leq \tilde{\varepsilon}$.
6:     **if** $L\|v\| \in [2\lambda_i - L\tilde{\varepsilon}, 2\lambda_i + L\tilde{\varepsilon}]$ **then**
7:         **return** $(\lambda_i, v, \emptyset)$.
8:     **else if** $L\|v\| > 2\lambda_i + L\tilde{\varepsilon}$ **then**
9:         **return** BinarySearch($\lambda_1 = \lambda_{i-1}, \lambda_2 = \lambda_i, \tilde{\varepsilon}$).
10:    **else if** $L\|v\| < 2\lambda_i - L\tilde{\varepsilon}$ **then**
11:        Let Power Method find vector $w$ that is 9/10-appx leading eigenvector of $(\mathbf{H} + \lambda_i\mathbf{I})^{-1}$:

$$\tfrac{9}{10}\lambda_{\max}((\mathbf{H} + \lambda_i\mathbf{I})^{-1}) \leq w^\top(\mathbf{H} + \lambda_i\mathbf{I})^{-1}w \leq \lambda_{\max}((\mathbf{H} + \lambda_i\mathbf{I})^{-1}) \ .$$

12:        Compute a vector $\tilde{w}$ such that $\|\tilde{w} - (\mathbf{H} + \lambda_i\mathbf{I})^{-1}w\| \leq \hat{\varepsilon} \triangleq \frac{1}{60B}$.
13:        $\Delta \leftarrow \frac{1}{2}\frac{1}{\tilde{w}^\top w - \hat{\varepsilon}}$ .
14:        **if** $\Delta > \frac{1}{2\kappa}$ **then**
15:            $\tilde{\lambda}_{i+1} \leftarrow \lambda_i - \frac{\Delta}{2}$.
16:            **if** $\tilde{\lambda}_{i+1} > 0$ **then** $\lambda_{i+1} \leftarrow \tilde{\lambda}_{i+1}$ **else** $\lambda_{i+1} \leftarrow 0$
17:        **else**
18:            Use AppxPCA to find any unit vector $v_{\min}$ such that $v_{\min}^\top \mathbf{H} v_{\min} \leq \lambda_{\min}(\mathbf{H}) + \frac{1}{10\kappa}$.
19:            Flip the sign of $v_{\min}$ so that $g^\top v_{\min} \leq 0$.
20:            **return** $(\lambda_i, v, v_{\min})$.
21:        **end if**
22:    **end if**
23: **end for**

---

**Algorithm 3** BinarySearch($\lambda_1, \lambda_2, \tilde{\varepsilon}$) (binary search subroutine)

**Input:** $\lambda_1 \geq \lambda_2$, $L\|(\mathbf{H}+\lambda_1\mathbf{I})^{-1}g\| \leq 2\lambda_1$, $L\|(\mathbf{H}+\lambda_2\mathbf{I})^{-1}g\| \geq 2\lambda_2$, $\lambda_2 + \lambda_{\min}(\mathbf{H}) > 0$
**Output:** $(\lambda, v, \emptyset)$
1: **for** $t = 1$ **to** $\infty$ **do**
2:     $\lambda_{\mathrm{mid}} \leftarrow \frac{\lambda_1+\lambda_2}{2}$
3:     Compute vector $v$ such that $\|v + (\mathbf{H}+\lambda_{\mathrm{mid}}\mathbf{I})^{-1}g\| \leq \tilde{\varepsilon}/2$
4:     **if** $L\|v\| \in [2\lambda_{\mathrm{mid}} - L\tilde{\varepsilon}, 2\lambda_{\mathrm{mid}} + L\tilde{\varepsilon}]$ **then**
5:         **return** $(\lambda_{\mathrm{mid}}, v, \emptyset)$
6:     **else if** $L\|v\| + L\tilde{\varepsilon} \leq 2\lambda_{\mathrm{mid}}$ **then**
7:         $\lambda_1 \leftarrow \lambda_{\mathrm{mid}}$
8:     **else if** $L\|v\| - L\tilde{\varepsilon} \geq 2\lambda_{\mathrm{mid}}$ **then**
9:         $\lambda_2 \leftarrow \lambda_{\mathrm{mid}}$
10:    **end if**
11: **end for**

*Correctness Proof of Main Lemma 3.* We carefully verify these sufficient conditions:

- Lemma 7.2 implies $\lambda_i \in [0, 2B]$.
- $\lambda_i + \lambda_{\min}(\mathbf{H}) \geq \frac{3}{10\kappa}$ from Lemma 7.2 implies $\|(\mathbf{H} + \lambda_i \mathbf{I})^{-1}\| \leq 4\kappa$. It is now immediate that the choice of $\tilde{\varepsilon}$ on Line 2 satisfies the Condition (6.1) in the assumption of Main Lemma 1.
- Since $\tilde{\varepsilon} \leq \frac{1}{10\kappa L}$ and $\lambda_i + \lambda_{\min}(\mathbf{H}) \geq \frac{3}{10\kappa}$ it follows that $(\mathbf{H} + (\lambda_i - L\tilde{\varepsilon})\mathbf{I})^{-1} \succ 0$ which proves Condition (6.2) in Main Lemma 1.
- We now verify Case 1 and 2 in the assumption of Main Lemma 1. At the beginning of the algorithm, our choice $\lambda_0 = 2B$ ensures (using Proposition 5.3) that $L\|(\mathbf{H} + \lambda_0 \mathbf{I})^{-1} g\| < 2\lambda_0$. Let us now consider the various places where the algorithm outputs:
    - If FastCubicMin terminates at Line 7, then we have $\|v + (\mathbf{H} + \lambda_i \mathbf{I})^{-1} g\| \leq \tilde{\varepsilon}$ and additionally
    $$L\|(\mathbf{H} + \lambda_i \mathbf{I})^{-1} g\| \in \big[L\|v\| - L\tilde{\varepsilon},\ L\|v\| + L\tilde{\varepsilon}\big] \subseteq [2\lambda_i - 2L\tilde{\varepsilon},\ 2\lambda_i + 2L\tilde{\varepsilon}] \ .$$
    Therefore, the output meets Case 1 requirement of Main Lemma 1 with $\lambda = \lambda_i$.
    - If FastCubicMin terminates at Line 9, then $L\|(\mathbf{H} + \lambda_i \mathbf{I})^{-1} g\| > L\|v\| - L\tilde{\varepsilon} \geq 2\lambda_i$ . Obviously, we must have $i \geq 1$ in this case because $L\|(\mathbf{H} + \lambda_0 \mathbf{I})^{-1} g\| < 2\lambda_0$. Therefore, Line 10 must have been reached at the previous iteration, so it implies $L\|(\mathbf{H} + \lambda_{i-1} \mathbf{I})^{-1} g\| < 2\lambda_{i-1}$ . Together, these two imply that we can call BinarySearch with $(\lambda_{i-1}, \lambda_i)$. Owing to Fact 7.1, the subroutine outputs a pair $(\lambda, v)$ satisfying the Case 1 requirement of Main Lemma 1.
    - If FastCubicMin terminates on Line 20, we verify that Case 2 of Main Lemma 1 with $\lambda = \lambda_i$ holds. We first have
    $$L\|(\mathbf{H} + \lambda_i \mathbf{I})^{-1} g\| \leq L\|v\| + L\tilde{\varepsilon} \leq 2\lambda_i \ .$$
    By Corollary 5.2, we also have that $\lambda_i \geq \lambda^*$. Lemma 7.2 tells us $\lambda_i$ satisfies $\lambda_i + \lambda_{\min}(\mathbf{H}) \leq \frac{1}{\kappa}$. Vector $v$ satisfies $\|v + (\mathbf{H} + \lambda_i \mathbf{I})^{-1} g\| \leq \tilde{\varepsilon}$. Vector $v_{\min}$ satisfies $v_{\min}^\top \mathbf{H} v_{\min} \leq \lambda_{\min}(\mathbf{H}) + \frac{1}{10\kappa}$ .

In sum, we have verified that all the assumptions of Main Lemma 1 hold. $\square$

**Final Proof of Theorem 2.** Theorem 2 is a direct corollary of our main lemmas. Main Lemma 3 ensures that the assumptions of Main Lemma 1 and Main Lemma 2 both hold. Now, using the special choice of $\kappa$ in FastCubic, Theorem 2-a immediately comes from Main Lemma 1; Theorem 2-b immediately comes from Main Lemma 2; and Theorem 2-c immediately comes from Main Lemma 3. This finishes the proof of Theorem 2. $\square$

## Acknowledgements

We thank Ben Recht for very helpful suggestions and corrections to a previous version. Z. Allen-Zhu is supported by an NSF Grant, no. CCF-1412958, and a Microsoft Research Grant, no. 0518584. Any opinions, findings and conclusions or recommendations expressed in this material are those of the author(s) and do not necessarily reflect the views of NSF or Microsoft.

# Appendix

## A  Computing Hessian-Vector Product in Linear Time

In this section we sketch the intuition regarding why Hessian-vector products can be computed in linear time in many interesting (especially machine learning) problems. We start by showing that the gradient can be computed in linear time. The algorithm is often referred to as back-propagation, which dates back to Werbos's PhD thesis [35], and has been popularized by Rumelharte *et al.* [31] for training neural networks.

**Claim A.1** (back-propagation, informally stated). *Suppose a real-valued function $f : \mathbb{R}^d \to \mathbb{R}$ can be evaluated by a differentiable circuit of size $N$. Then, the gradient $\nabla f$ can be computed in time $O(N + d)$ (using a circuit of size $O(N + d)$).* [4]

The claim follows from simple induction and chain-rule, and is left to the readers. In the training of neural networks, often the size of circuits that computes the objective $f$ is proportional to (or equal to) the number of parameters $d$. Thus the gradient $\nabla f$ can be computed in time $O(d)$ using a circuit of size $d$.

Next, we consider computing $\nabla^2 f(x) \cdot v$ where $v \in \mathbb{R}^d$. Let $g(x) := \langle \nabla f(x), v \rangle$ be a function from $\mathbb{R}^d$ to $\mathbb{R}$. Then, we see that if suffices to compute the gradient of $g$, since
$$\nabla^2 f(x) \cdot v = \nabla g(x).$$
We observe that $g(x)$ can be evaluated in linear time using circuit of size $O(d)$ since we've shown $\nabla f(x)$ can. Thus, using Claim A.1 again on function $g$, [5] we conclude that $\nabla g(x)$ can also be computed in linear time.

## B  Proof of Lemma B.1 and Corollary 4.1

**Lemma B.1.** *For all $h' \in \mathbb{R}^d$, it satisfies*
$$\|\nabla f(x_t+h')\| \leq L\|h'\|^2 + \|\nabla m_t(h')\| \quad \text{and} \quad \lambda_{\min}(\nabla^2 f(x_t+h')) \geq -\left(\frac{3L^2 \max\{0, -m_t(h^*)\}}{2}\right)^{1/3} - L\|h'\|.$$

*Proof of Lemma B.1.* Let us denote by $g = \nabla f(x_t)$ and $\mathbf{H} = \nabla^2 f(x_t)$ in this proof. We begin by proving the first order condition. Note that we have
$$\nabla m_t(h) = g + \mathbf{H}h + \tfrac{L}{2}\|h\|h.$$
Recall $h^*$ is a minimizer of $\operatorname{argmin} m_t(h)$. The characterization result in Lemma 5.1 shows $\mathbf{H} + \frac{L\|h^*\|}{2}\mathbf{I} \succeq 0$, and thus
$$g^\top h^* + (h^*)^\top \mathbf{H} h^* + \frac{L}{2}\|h^*\|^3 = \nabla m_t(h^*)^\top h^* = 0 \tag{B.1}$$
$$(h^*)^\top \mathbf{H} h^* + \frac{L}{2}\|h\|^3 = (h^*)^\top \left(\mathbf{H} + \frac{L\|h^*\|}{2}\mathbf{I}\right) h^* \geq 0. \tag{B.2}$$
They imply
$$m_t(h^*) = g^\top h^* + \frac{(h^*)^\top \mathbf{H} h^*}{2} + \frac{L\|h^*\|^3}{6} \stackrel{①}{=} -\frac{(h^*)^\top \mathbf{H} h^*}{2} - \frac{L}{3}\|h^*\|^3$$
$$\stackrel{②}{\leq} \frac{L}{4}\|h^*\|^3 - \frac{L}{3}\|h^*\|^3 = -\frac{L}{12}\|h^*\|^3 \tag{B.3}$$
where ① uses (B.1) and ② uses (B.2).

We compute the norm of the gradient at a point $x_t + h'$ for any $h' \in \mathbb{R}^d$:
$$\|\nabla f(x_t + h')\| \leq \|\nabla f(x_t + h') - \nabla m_t(h')\| + \|\nabla m_t(h')\|$$
$$= \left\|\nabla f(x_t) + \int_0^1 \nabla^2 f(x_t + \tau h')h' d\tau - \left(g + \mathbf{H}h' + \frac{L}{2}\|h'\|h'\right)\right\| + \|\nabla m_t(h')\|$$

---
[4]Technically, we assume that the gradient of each gate can be computed in $O(1)$ time
[5]We assume here that the original circuits are twice differentiable

$$
\leq \left\| \int_0^1 \left( \nabla^2 f(x_t + \tau h') - \mathbf{H} \right) h' d\tau \right\| + \frac{L}{2} \|h'\|^2 + \|\nabla m_t(h')\|
$$
$$
\stackrel{③}{\leq} L\|h'\|^2 \int_0^1 \tau d\tau + \frac{L}{2}\|h'\|^2 + \|\nabla m_t(h')\| = L\|h'\|^2 + \|\nabla m_t(h')\| \tag{B.4}
$$
where ③ follows from the Lipschitz continuity on the Hessian (2.1). This proves the first conclusion of the lemma.

As for the second-order condition, we first note that for all $h' \in \mathbb{R}^d$, by the Lipschitz continuity on the Hessian (2.1), we have $\|\nabla^2 f(x_t + h') - \nabla^2 f(x_t)\| \leq L\|h'\|$. However, this implies
$$
\lambda_{\min}(\nabla^2 f(x_t + h')) \geq \lambda_{\min}(\nabla^2 f(x_t)) - L\|h'\| \ . \tag{B.5}
$$
because if two matrices $\mathbf{A}$ and $\mathbf{B}$ satisfies $\|\mathbf{A}-\mathbf{B}\| \leq p$, then it must satisfy $\left|\lambda_{\min}(\mathbf{A}) - \lambda_{\min}(\mathbf{B})\right| \leq p$ as well. We consider two cases: if $\lambda_{\min}(\nabla^2 f(x_t)) \geq 0$, then we have
$$
\lambda_{\min}(\nabla^2 f(x_t + h')) \geq -L\|h'\| \ . \tag{B.6}
$$
Otherwise, we consider the case where $\lambda_{\min}(\nabla^2 f(x_t)) = \lambda_{\min}(\mathbf{H}) < 0$. Let $\nu_d$ be the normalized eigenvector corresponding to $\lambda_{\min}(\mathbf{H})$, and define
$$
\tilde{h} \triangleq \mathrm{sign}(g^\top \nu_d) \cdot \frac{2\lambda_{\min}(\mathbf{H})}{L} \nu_d \ .
$$
We calculate $m_t(\tilde{h})$ as follows:
$$
m_t(\tilde{h}) = g^\top \tilde{h} + \frac{\tilde{h}^\top \mathbf{H} \tilde{h}}{2} + \frac{L}{6}\|\tilde{h}\|^3 \leq \frac{\tilde{h}^\top \mathbf{H} \tilde{h}}{2} + \frac{L}{6}\|\tilde{h}\|^3 = \frac{2(\lambda_{\min}(\mathbf{H}))^2}{L^2}\nu_d^\top \mathbf{H} \nu_d + \frac{4|\lambda_{\min}(\mathbf{H})|^3}{3L^2}
$$
$$
\stackrel{①}{=} \frac{2(\lambda_{\min}(\mathbf{H}))^3}{L^2} + \frac{4|\lambda_{\min}(\mathbf{H})|^3}{3L^2} \stackrel{②}{=} \frac{2(\lambda_{\min}(\mathbf{H}))^3}{3L^2} \ , \tag{B.7}
$$
where ① uses $\nu_d^\top \mathbf{H} \nu_d = \lambda_{\min}(\mathbf{H}) < 0$, and ② uses the assumption that $\lambda_{\min}(\mathbf{H}) < 0$. Since by definition $m_t(h^*) \leq m_t(\tilde{h})$, we can deduce from inequality (B.7) that
$$
\lambda_{\min}(\nabla^2 f(x_t)) = \lambda_{\min}(\mathbf{H}) \geq -\left(\frac{3L^2|m_t(h^*)|}{2}\right)^{1/3} \ . \tag{B.8}
$$
Now we put together inequalities (B.5) and (B.8), and obtain
$$
\lambda_{\min}(\nabla^2 f(x_t + h')) \geq -\left(\frac{3L^2|m_t(h^*)|}{2}\right)^{1/3} - L\|h'\| \ . \tag{B.9}
$$
Combining (B.6) and (B.9) we finish the proof of Lemma B.1. $\square$

**Corollary 4.1.** *If $m_t(h^*) \geq -\frac{\varepsilon^{3/2}}{800\sqrt{L}}$ and $h'$ is an approximate minimizer of $m_t(h)$ satisfying*
$$
\|h'\| \leq \|h^*\| + \tfrac{\sqrt{\varepsilon}}{4\sqrt{L}} \quad \text{and} \quad \|\nabla m_t(h')\| \leq \tfrac{\varepsilon}{2} \ ,
$$
*then we have that $\|\nabla f(x_t + h')\| \leq \varepsilon$ and $\lambda_{\min}(\nabla^2 f(x_t + h')) \geq -\sqrt{L\varepsilon}$.*

*Proof of Corollary 4.1.* First of all, our assumption that $m_t(h^*) \geq -\frac{\varepsilon^{3/2}}{800\sqrt{L}}$, along with inequality (B.3), tells us that $\|h^*\| \leq \frac{\sqrt{\varepsilon}}{4\sqrt{L}}$. This, together with our assumption on $\|h'\|$, implies $\|h'\| \leq \frac{\sqrt{\varepsilon}}{2\sqrt{L}}$. Since we also assume $\|\nabla m_t(h')\| \leq \frac{\varepsilon}{2}$, we have from Lemma B.1 that
$$
\|\nabla f(x_t + h')\| \leq L\|h'\|^2 + \|\nabla m_t(h')\| \leq \frac{\varepsilon}{4} + \frac{\varepsilon}{2} \leq \varepsilon \ .
$$
For the second-order condition, we can again apply Lemma B.1 to get
$$
\lambda_{\min}(\nabla^2 f(x_t + h')) \geq -\left(\frac{3L^2 \max\{0, -m_t(h^*)\}}{2}\right)^{1/3} - L\|h'\| \geq -\left(\frac{3L^{3/2}\varepsilon^{3/2}}{1600}\right)^{1/3} - \frac{\sqrt{L\varepsilon}}{2} \geq -\sqrt{L\varepsilon} \ .
$$

# C Proof of Lemma 5.1 and Corollary 5.2

We begin by proving a few lemmas that characterize the system of equations.

**Lemma C.1.** *Consider the following system of equations/inequalities in variables $\lambda, h$:*

$$\mathbf{H} + \lambda \mathbf{I} \succeq 0 \ , \qquad (\mathbf{H} + \lambda \mathbf{I})h = -g \ , \qquad \|h\| = \frac{2\lambda}{L} \ . \tag{C.1}$$

*The following statements hold for any solution $(\lambda', h')$ of the above system:*

- *There is a unique value $\lambda'$ that satisfies the above equations. $\lambda'$ is such that $\lambda' \geq -\lambda_{\min}(\mathbf{H})$.*
- *If $\lambda' > -\lambda_{\min}(\mathbf{H})$, then the corresponding $h'$ is also unique and is given by $h' = -(\mathbf{H}+\lambda\mathbf{I})^{-1}g$.*
- *If $\lambda' = -\lambda_{\min}(\mathbf{H})$ then $g^\top v = 0$ for any vector $v$ belonging to the eigenspace corresponding to $\lambda_{\min}(\mathbf{H})$. Subsequently we also have that the corresponding $h'$ is of the form*

$$h' = -(\mathbf{H} + \lambda I)^+ g + \gamma v$$

*for some $\gamma$ and $v$ in the lowest eigenspace of $\mathbf{H}$.*

*Proof of Lemma C.1.* Note that $\mathbf{H} + \lambda \mathbf{I} \succeq 0$ ensures that for any solution $\lambda'$, we have $\lambda' \geq -\lambda_{\min}(\mathbf{H})$. Furthermore, for any $\lambda' > -\lambda_{\min}(\mathbf{H})$, the corresponding $h$ is uniquely defined by $h = (\mathbf{H} + \lambda \mathbf{I})^{-1}g$ since $\mathbf{H} + \lambda' \mathbf{I}$ is invertible. If indeed $\lambda' = -\lambda_{\min}(\mathbf{H})$, then we have that the equation $(\mathbf{H} - \lambda_{\min}(\mathbf{H})\mathbf{I})h = -g$ has a solution. This implies that $g$ has no component in the null space of $\mathbf{H} - \lambda_{\min}(\mathbf{H})\mathbf{I}$, or equivalently that it has no component in the eigenspace corresponding to $\lambda_{\min}(\mathbf{H})$. We also have that every solution of $(\mathbf{H} - \lambda_{\min}(\mathbf{H})\mathbf{I})h = -g$ is necessarily of the form

$$h = -(H - \lambda_{\min}\mathbf{I})^+ g + \gamma v$$

for some $\gamma$ and $v$ in the lowest eigenspace of $\mathbf{H}$.

We will now prove the uniqueness of $\lambda'$ by contradiction. Consider two distinct values of $\lambda_1, \lambda_2$ that satisfy the system (C.1). If both $\lambda_1, \lambda_2 > -\lambda_{\min}(\mathbf{H})$ we get that

$$\|(\mathbf{H} + \lambda_1 \mathbf{I})^{-1}g\| = \frac{2\lambda_1}{L} \quad \text{and} \quad \|(\mathbf{H} + \lambda_2 \mathbf{I})^{-1}g\| = \frac{2\lambda_2}{L} \ .$$

Now note that $\|(\mathbf{H} + \lambda \mathbf{I})^{-1}g\|$ is a strictly decreasing function over the domain $\lambda \in (-\lambda_{\min}(H), \infty)$ and $\frac{2\lambda}{L}$ is strictly increasing over the same domain. Therefore the above two equations cannot be satisfied for two distinct $\lambda_1, \lambda_2 > -\lambda_{\min}(\mathbf{H})$ which is a contradiction. Suppose now without loss of generality that $\lambda_1 = -\lambda_{\min}(\mathbf{H})$. Then we have that the corresponding solution is of the form

$$h = -(H + \lambda\mathbf{I})^+ g + \gamma v$$

for some $\gamma$ and $v$ in the lowest eigenspace of $\mathbf{H}$ and $g$ has no component in the lowest eigenspace of $\mathbf{H}$. It follows that $\|(\mathbf{H} - \lambda_{\min}(\mathbf{H})\mathbf{I})^+ g\| \geq \|(\mathbf{H} + \lambda\mathbf{I})^{-1}g\|$ for any $\lambda > -\lambda_{\min}(\mathbf{H})$. By a similar argument as in the first case, we can now see that the following conditions,

$$\|(\mathbf{H} + \lambda_1 \mathbf{I})^+ g + \gamma v_{\min}(\mathbf{H})\| = \frac{2\lambda_1}{L} \quad \text{and} \quad \|(\mathbf{H} + \lambda_2 \mathbf{I})^{-1}g\| = \frac{2\lambda_2}{L} \ ,$$

cannot both be satisfied for $\lambda_2 > \lambda_1 = -\lambda_{\min}(\mathbf{H})$, giving us a contradiction. This finishes the proof of Lemma C.1. □

**Lemma C.2.** *Let $(\lambda, h)$ be a solution of the system (C.1). Then we have that*

$$m(h) = -\frac{1}{2}g^\top (\mathbf{H} + \lambda \mathbf{I})^+ g - \frac{2\lambda^3}{3L^2} \ .$$

*Proof of Lemma C.2.* By the definition of the system (C.1), any solution $\lambda, h$ to the system should be such that there exists some $\gamma$ such that
$$h = (H + \lambda \mathbf{I})^+ g + \gamma v_0$$
where $v_0$ is in the null space of $\mathbf{H} + \lambda \mathbf{I}$ if it exists; otherwise $\gamma = 0$. This gives us the following:
$$\begin{aligned} m(h) &= g^\top h + \frac{h^\top \mathbf{H} h}{2} + \frac{L}{6}\|h\|^3 \\ &\overset{\text{①}}{=} -\frac{1}{2}h^\top(\mathbf{H} + \lambda\mathbf{I})h - \frac{\lambda}{2}\|h\|^2 + \frac{L}{6}\|h\|^3 \\ &\overset{\text{②}}{=} -\frac{1}{2}g^\top(\mathbf{H} + \lambda\mathbf{I})^+ g - \frac{2\lambda^3}{3L^2} \ . \end{aligned}$$
Equality ① follows because $(\mathbf{H} + \lambda\mathbf{I})h = -g$. Equality ② follows because $h = (\mathbf{H} + \lambda\mathbf{I})^+ g + \gamma v_0$ and $\|h\| = \frac{2\lambda}{L}$. □

**Lemma 5.1.** $h^*$ *is a minimizer of $m(h)$ if and only if there exists $\lambda^* \geq 0$ such that*
$$\mathbf{H} + \lambda^*\mathbf{I} \succeq 0 \ , \qquad (\mathbf{H} + \lambda^*\mathbf{I})h^* = -g \ , \qquad \|h^*\| = \frac{2\lambda^*}{L} \ .$$
*The objective value in this case is given by*
$$m(h^*) = -\frac{1}{2}g^\top(\mathbf{H} + \lambda^*\mathbf{I})^+ g - \frac{2(\lambda^*)^3}{3L^2} \leq 0 \ .$$

*Proof of Lemma 5.1.* We first compute that
$$\nabla m(h) = g + \mathbf{H}h + \frac{L}{2}\|h\|h \quad \text{and} \quad \nabla^2 m(h) = \mathbf{H} + \frac{L}{2}\|h\|\mathbf{I} + \frac{L}{2}\|h\|\left(\frac{h}{\|h\|}\right)\left(\frac{h}{\|h\|}\right)^\top \ .$$
For the forward direction, suppose $h^*$ is a minimizer of $m(h)$. Let $\lambda^* = \frac{L}{2}\|h^*\|$. Then, the necessary conditions $\nabla m(h^*) = 0$ and $\nabla^2 m(h^*) \succeq 0$ can be written as
$$g + (\mathbf{H} + \lambda^*\mathbf{I})h^* = 0 \quad \text{and} \quad w^\top \left(\mathbf{H} + \lambda^*\mathbf{I} + \lambda^*\left(\frac{h^*}{\|h^*\|}\right)\left(\frac{h^*}{\|h^*\|}\right)^\top\right) w \geq 0, \forall w \in \mathbb{R}^n. \quad (\text{C.2})$$
From this we see $(\mathbf{H} + \lambda^*\mathbf{I})h^* = -g$ and $\|h^*\| = \frac{2\lambda^*}{L}$, and the only thing left to verify is $\mathbf{H} + \lambda^*\mathbf{I} \succeq 0$.

Note that if $h^* = 0$, then the second inquality in (C.2) directly implies $\mathbf{H} + \lambda^*\mathbf{I} \succeq 0$. Thus, we only need to focus on $h^* \neq 0$. We want to show that $w^\top(\mathbf{H} + \lambda^*\mathbf{I})w \geq 0$ for every $w \in \mathbb{R}^d$. Now, if $w^\top h^* = 0$ then this trivially follows from (C.2), so it suffices to focus on those $w$ that satisfies $w^\top h^* \neq 0$.

Since $w$ and $h^*$ are not orthogonal, there exists $\gamma \in \mathbb{R}\setminus\{0\}$ such that $\|h^* + \gamma w\| = \|h^*\|$. (This can be done by squaring both sides and solving the linear system in $\lambda$.) Squaring both sides we have
$$(\gamma w)^\top h^* + \frac{\gamma^2 \|w\|^2}{2} = 0 \ . \tag{C.3}$$
Now we bound the difference
$$\begin{aligned} m(h^* + \gamma w) - m(h^*) &= g^\top((h^* + \gamma w) - h^*) + \frac{(h^* + \gamma w)^\top \mathbf{H}(h^* + \gamma w)}{2} - \frac{h^{*}\mathbf{H}h^*}{2} \\ &\overset{\text{①}}{=} (h^* - (h^* + \gamma w))^\top(\mathbf{H} + \lambda^*\mathbf{I})h^* + \frac{(h^* + \gamma w)^\top \mathbf{H}(h^* + \gamma w)}{2} - \frac{h^{*}\mathbf{H}h^*}{2} \\ &\overset{\text{②}}{=} \frac{\lambda^* \gamma^2}{2}\|w\|^2 + (h^* - (h^* + \gamma w))^\top \mathbf{H}h^* + \frac{(h^* + \gamma w)^\top \mathbf{H}(h^* + \gamma w)}{2} - \frac{h^{*}\mathbf{H}h^*}{2} \\ &= \frac{\lambda^* \gamma^2}{2}\|w\|^2 + \frac{h^{*}\mathbf{H}h^*}{2} - (h^* + \gamma w)^\top \mathbf{H}h^* + \frac{(h^* + \gamma w)^\top \mathbf{H}(h^* + \gamma w)}{2} \end{aligned}$$

$$= \frac{\lambda^* \gamma^2}{2}\|w\|^2 + \frac{\gamma^2}{2}w^\top \mathbf{H} w = \frac{\gamma^2}{2}w^\top(\mathbf{H}+\lambda^*\mathbf{I})w \ , \tag{C.4}$$

where ① and ② follow from (C.2) and (C.3), respectively. Since $h^*$ is a minimizer of $m(h)$, we immediately have

$$m(h^*+\gamma w) - m(h^*) = \frac{\gamma^2}{2}w^\top(\mathbf{H}+\lambda^*\mathbf{I})w \geq 0,$$

and we conclude that $(\mathbf{H}+\lambda^*\mathbf{I}) \succeq 0$.

For the backward direction, we will make use Lemma C.1 and Lemma C.2. First we note that the function $m(h)$ is continuous and bounded from below, and there exists at least one minimizer $h^*$. Suppose now there exists a $\lambda^*$ and a corresponding $h^*$ such that $(\lambda^*, h^*)$ is a solution to the system C.1. The backward direction requires us to prove that $h^*$ must be a minimizer of $m(h)$. By Lemma C.1 we get the following two cases.

We prove the backward direction by showing that the conditions in Equation C.2 determine the minimizer up to its norm. To this end we will use Lemma C.1 and Lemma C.2.

First we note that the function $m(h)$ is continuous, bounded from below, and tends to $+\infty$ when $\|h\| \to \infty$, so there exists at least one minimizer $h^*$.

Suppose now there exists a $\lambda^*$ and a corresponding $h^*$ such that $(\lambda^*, h^*)$ is a solution to the system (C.1). The backward direction requires us to prove that $h^*$ must be a minimizer of $m(h)$. By Lemma C.1 we get the following two cases.

- If $\lambda^* > -\lambda_{\min}(\mathbf{H})$ then $(\lambda^*, h^*)$ is the only solution to the system (C.1). By the proof of the forward direction we see that any minimizer of $m(h)$ must satisfy system (C.1) and therefore $h^*$ must be the minimizer.

- If above is not the case, then $\lambda^* = -\lambda_{\min}(\mathbf{H})$. Let $h'$ be any minimizer of $m(h)$. Lemma C.1 and the proof of the forward direction ensures that $(\lambda^*, h')$ also satisfies the system (C.1). By Lemma C.2 we get $m(h^*) = m(h')$ and therefore $h^*$ is a minimizer too. □

**Corollary 5.2.** *This value $\lambda^*$ is unique, and for every $\lambda$ satisfying $\mathbf{H}+\lambda\mathbf{I} \succ 0$, we have*

$$\|(\mathbf{H}+\lambda\mathbf{I})^{-1}g\| > \frac{2\lambda}{L} \iff \lambda^* > \lambda \quad \text{and} \quad \|(\mathbf{H}+\lambda\mathbf{I})^{-1}g\| < \frac{2\lambda}{L} \iff \lambda^* < \lambda \ .$$

*Proof of Corollary 5.2.* The uniqueness of $\lambda^*$ follows from Lemma C.1. To prove the second part we first make some observations about the function

$$p(y) \triangleq \frac{2y}{L} - \|(\mathbf{H}+y\mathbf{I})^{-1}g\|$$

defined on the domain $y \in (-\lambda_{\min}(H), \infty)$. Note that $p(y)$ is continuous and strictly increasing over the domain and $p(y) \to \infty$ as $y \to \infty$.

The corollary requires us to show that

$$p(\lambda) < 0 \iff \lambda^* > \lambda \quad \text{and} \quad p(\lambda) > 0 \iff \lambda^* < \lambda \ .$$

We begin by showing the first equivalence. To see the backward direction note that if $\lambda^* > \lambda > -\lambda_{\min}(\mathbf{H})$, by the characterization of $\lambda^*$ in Lemma C.1 we have that $\|(\mathbf{H}+\lambda^*\mathbf{I})^{-1}g\| = \frac{2\lambda^*}{L}$ i.e. $p(\lambda^*) = 0$ which implies that $p(\lambda) < 0$ as $p(y)$ is a strictly increasing function. For the forward direction note that since $p(y)$ is continuous and strictly increasing we see that the range of the function contains $[p(\lambda), \infty)$. Since $p(\lambda) < 0$ there must exist a $\lambda^* > \lambda$ such that $p(\lambda^*) = 0$ which by the characterization in Lemma C.1 finishes the proof.

Now we will prove that $p(\lambda) > 0 \iff \lambda^* < \lambda$. To see the forward direction note that if $\lambda^* \geq \lambda$ then $p(\lambda^*) = 0$ and $p(\lambda) > 0$ which contradicts the fact that $p(y)$ is strictly increasing. For the backward direction we consider two cases. Firstly if $\lambda^* > -\lambda_{\min}(\mathbf{H})$ the conclusion follows similarly by the monotonicity of $p(y)$. If $\lambda^* = -\lambda_{\min}$ then by Lemma C.1, we have that $g$ has no component in the lowest eigenspace of $\mathbf{H}$ and therefore if we extend $p(y)$ to $-\lambda_{\min}(\mathbf{H})$ by defining

$$p(-\lambda_{\min}(\mathbf{H})) \triangleq \frac{-2\lambda_{\min}(\mathbf{H})}{L} - \|(\mathbf{H} - \lambda_{\min}(\mathbf{H})\mathbf{I})^+ g\|$$

we get that $p(y)$ is increasing in the domain $y \in [-\lambda_{\min}(\mathbf{H}), \infty)$. Now from the characterization of the solution in Lemma C.1 we can see that $p(-\lambda_{\min}(\mathbf{H})) \geq 0$ and therefore by monotonicity $p(\lambda) > 0$. This finishes the proof. □

# D  Proof of Main Lemma 1

## D.1  Proof of Claim 6.1

**Claim 6.1.** *If $\lambda$ and $v$ satisfy Case 1 and $\tilde{\varepsilon}$ satisfies (6.1), then $m(v) \leq m(h^*) + \frac{1}{250\kappa^3 L^2}$*

*Proof of Claim 6.1.* Note that by the conditions of the theorem we have that $(\mathbf{H}+(\lambda-L\tilde{\varepsilon})\mathbf{I})^{-1} \succ eq$ and

$$L\|(\mathbf{H} + (\lambda - L\tilde{\varepsilon})\mathbf{I})^{-1}g\| \geq 2\lambda - 2L\tilde{\varepsilon} \quad \text{and} \quad L\|(\mathbf{H} + (\lambda + L\tilde{\varepsilon})\mathbf{I})^{-1}g\| \leq 2\lambda - 2L\tilde{\varepsilon} ,$$

according to Corollary 5.2 we must have

$$\lambda^* \in [\lambda - L\tilde{\varepsilon}, \lambda + L\tilde{\varepsilon}] \tag{D.1}$$

This also implies (using our assumption on $\tilde{\varepsilon}$)

$$L\|v\| \leq [2\lambda^* - 5L\tilde{\varepsilon}, 2\lambda^* + 5L\tilde{\varepsilon}] .$$

Next, consider the value $m(v)$

$$m(v) = g^\top v + \frac{v^\top H v}{2} + \frac{L}{6}\|v\|^3 = g^\top v + \frac{v^\top (\mathbf{H} + \lambda \mathbf{I})v}{2} - \|v\|^2 \left(\frac{\lambda}{2} - \frac{L\|v\|}{6}\right) . \tag{D.2}$$

We bound the two parts on the right hand side of (D.2) separately. The first part

$$g^\top v + \frac{v^\top (\mathbf{H} + \lambda \mathbf{I})v}{2} \leq -\frac{g^\top (\mathbf{H} + \lambda \mathbf{I})^{-1}g}{2} + \|g\|\tilde{\varepsilon} + \|(\mathbf{H} + \lambda \mathbf{I})^{-1}g\|\tilde{\varepsilon} + \frac{\|(\mathbf{H} + \lambda \mathbf{I})^{-1}\|\tilde{\varepsilon}^2}{2}$$

$$\overset{①}{\leq} -\frac{g^\top (\mathbf{H} + \lambda \mathbf{I})^{-1}g}{2} + \frac{1}{1000\kappa^3 L^2} \tag{D.3}$$

$$\overset{②}{\leq} -\frac{g^\top (\mathbf{H} + \lambda^* \mathbf{I})^{-1}g}{2} + L\|g\|^2 \|(\mathbf{H} + \lambda \mathbf{I})^{-1}\| \|(\mathbf{H} + (\lambda + 2L\tilde{\varepsilon})\mathbf{I})^{-1}\|\tilde{\varepsilon} + \frac{1}{1000\kappa^3 L^2}$$

$$\overset{③}{\leq} -\frac{g^\top (\mathbf{H} + \lambda^* \mathbf{I})^{-1}g}{2} + \frac{1}{500\kappa^3 L^2}$$

Above, inequalities ① and ③ use the assumption on $\tilde{\varepsilon}$ in (6.1), and inequality ② uses

$$-(\mathbf{H} + \lambda \mathbf{I})^{-1} \preceq -(\mathbf{H} + (\lambda^* + L\tilde{\varepsilon})\mathbf{I})^{-1}$$
$$= -(\mathbf{H} + \lambda^* \mathbf{I})^{-1} - L\tilde{\varepsilon}(\mathbf{H} + \lambda^* \mathbf{I})^{-1}(\mathbf{H} + (\lambda^* + L\tilde{\varepsilon})\mathbf{I})^{-1}$$

Note that $(\mathbf{H}+\lambda^*\mathbf{I})^{-1} \succ 0$ by Equations (D.1) and (6.2). The second part of (D.2) can be bounded as follows

$$\|v\|^2 \left(\frac{\lambda}{2} - \frac{L\|v\|}{6}\right) \geq \frac{(2\lambda^* - 5L\tilde{\varepsilon})^2}{L^2} \left(\frac{\lambda^* - L\tilde{\varepsilon}}{2} - \frac{2\lambda^* + 5L\tilde{\varepsilon}}{6}\right)$$

$$\geq \frac{2(\lambda^*)^3}{3L^2} - 1000\tilde{\varepsilon}L(\lambda^*)^2 \overset{①}{\geq} \frac{2(\lambda^*)^3}{3L^2} - \frac{1}{500\kappa^3 L^2}$$

Above, inequality ① uses $\lambda^* \leq B$ (owing to Proposition 5.3) and our assumption on $\tilde{\varepsilon}$ from (6.1). Putting these together we get that
$$m(v) \leq m(h^*) + \frac{1}{250\kappa^3 L^2} \quad . \qquad \square$$

## D.2 Proofs for Claims 6.2 and 6.3

For notational simplicity, let us rotate the space into the basis in the eigenspace of $\mathbf{H}$; let the $i$-th dimension correspond to the $i$-th largest eigenvalue $\lambda_i$ of $\mathbf{H}$. We have $\lambda_1 \geq \lambda_2 \ldots \geq \lambda_d = \lambda_{\min}$. Let $g_i$ denote the $i$-th coordinate of $g$ in this basis.

Lemma 5.1 implies
$$m(h^*) = -\frac{1}{2}\sum_i \frac{g_i^2}{\lambda_i + \lambda^*} - \frac{2(\lambda^*)^3}{3L^2} =: S_1 + S_2 - \frac{2(\lambda^*)^3}{3L^2} \quad . \tag{D.4}$$

where we denote by
$$S_1 = -\sum_{i:\lambda_i+\lambda^* \geq \frac{1}{\kappa}} \frac{g_i^2}{\lambda_i + \lambda^*} \qquad S_2 = -\sum_{i:0<\lambda_i+\lambda^* \leq \frac{1}{\kappa}} \frac{g_i^2}{\lambda_i + \lambda^*}$$

From Corollary 5.2 we can also obtain
$$\sum_{i:\lambda_i+\lambda^*>0} \frac{g_i^2}{(\lambda_i+\lambda^*)^2} \leq \frac{4(\lambda^*)^2}{L^2} \quad . \tag{D.5}$$

Now the assumption $\|(\mathbf{H} + \lambda\mathbf{I})^{-1}g\| \leq \frac{2\lambda}{L}$ is equivalent to
$$\sum_i \frac{g_i^2}{(\lambda_i + \lambda)^2} \leq \frac{4\lambda^2}{L^2} \tag{D.6}$$

We begin with a few auxiliary claims.

**Claim D.1.** *If $\lambda_{\min}(\mathbf{H}) \leq -\frac{1}{\kappa}$ then $S_2 \geq 1000 \cdot m\left(\frac{\lambda v_{\min}}{2L}\right)$*

*Proof of Claim D.1.* We compute that
$$S_2 = -\sum_{i:0<\lambda_i+\lambda^* \leq \frac{1}{\kappa}} \frac{g_i^2}{\lambda_i + \lambda^*} = -\sum_{i:0<\lambda_i+\lambda^* \leq \frac{1}{\kappa}} \frac{g_i^2(\lambda_i+\lambda^*)}{(\lambda_i+\lambda^*)^2} \geq -\frac{1}{\kappa}\sum_{i:0<\lambda_i+\lambda^* \leq \frac{1}{\kappa}} \frac{g_i^2}{(\lambda_i+\lambda^*)^2}$$
$$\overset{①}{\geq} -\frac{4}{\kappa L^2}(\lambda^*)^2 \overset{②}{\geq} -16\frac{|\lambda_{\min}|^3}{L^2} \quad . \tag{D.7}$$

Above, ① uses (D.5), and ② follows because we have $\lambda_{\min}(\mathbf{H}) \leq -\frac{1}{\kappa}$ in the assumption and have $\lambda^* \leq -\lambda_{\min}(\mathbf{H}) + \frac{1}{\kappa}$ in the assumption of Case 2 of Main Lemma 1.

Let us now consider the value of the vector $\frac{\lambda v_{\min}}{2L}$. We have that
$$m\left(\frac{\lambda v_{\min}}{2L}\right) = \frac{\lambda g^\top v_{\min}}{2L} + \frac{\lambda^2 v_{\min}^\top \mathbf{H} v_{\min}}{8L^2} + \frac{\lambda^3}{48L^2} \overset{①}{\leq} \frac{\lambda g^\top v_{\min}}{2L} + \frac{\lambda^2 \lambda_{\min}}{16L^2} + \frac{\lambda^3}{48L^2}$$
$$\overset{②}{\leq} \frac{\lambda g^\top v_{\min}}{2L} + \frac{\lambda^2 \lambda_{\min}}{16L^2} - \frac{\lambda^2 \lambda_{\min}}{24L^2} \leq \frac{\lambda g^\top v_{\min}}{2L} + \frac{\lambda^2 \lambda_{\min}}{48L^2}$$

Above, ① is because our assumption $\lambda_{\min}(\mathbf{H}) \leq -\frac{1}{\kappa}$ and assumption $v_{\min}\mathbf{H}v_{\min} \leq \lambda_{\min}(\mathbf{H}) + \frac{1}{10\kappa}$ together imply $v_{\min}\mathbf{H}v_{\min} \leq \frac{\lambda_{\min}}{2}$. ② follows from $\lambda_{\min}(\mathbf{H}) \leq -\frac{1}{\kappa}$ and $\lambda \leq -\lambda_{\min}(\mathbf{H}) + \frac{1}{\kappa}$.

Now, recall that the sign of $v_{\min}$ is chosen so $g^\top v_{\min}$ is non-positive, and therefore by our assumptions $\lambda_{\min}(\mathbf{H}) \leq -\frac{1}{\kappa}$ and $\lambda \leq -\lambda_{\min}(\mathbf{H}) + \frac{1}{\kappa}$, we get the following inequality:
$$m\left(\frac{\lambda v_{\min}}{2L}\right) \leq -\frac{|\lambda_{\min}|^3}{48L^2} \tag{D.8}$$
Putting inequalities (D.8) and (D.7) together finishes the proof of Claim D.1. □

We also show the following lemma, the proof of which can be seen from inequality (D.3), as part of the proof of Claim 6.1 above.

**Lemma D.2.** *If we have $\lambda, v$ such that*
$$L\|(\mathbf{H} + \lambda \mathbf{I})^{-1} g\| \leq 2\lambda \quad \text{and} \quad \|v + (\mathbf{H} + \lambda \mathbf{I})^{-1} g\| \leq \tilde{\varepsilon}$$
*with $\tilde{\varepsilon}$ satisfying condition (6.1) then we have that*
$$g^\top v + \frac{v^\top (\mathbf{H} + \lambda \mathbf{I}) v}{2} \leq -\frac{g^\top (\mathbf{H} + \lambda \mathbf{I})^{-1} g}{2} + \frac{1}{1000\kappa^3 L^2}$$

**Claim D.3.** $S_1 \geq 4m(v) - \frac{1}{250\kappa^3 L^2}$

*Proof of Claim D.3.* We have that
$$\begin{aligned}
m(v) &= g^\top v + \frac{v^\top(\mathbf{H}+\lambda\mathbf{I})v}{2} - \frac{\lambda}{2}\|v\|^2 + \frac{L}{6}\|v\|^3 \\
&\overset{①}{=} -\frac{g^\top(\mathbf{H}+\lambda\mathbf{I})^{-1}g}{2} - \|v\|^2\left(\frac{\lambda}{2} - \frac{L}{6}\|v\|\right) + \frac{1}{1000\kappa^3 L^2} \\
&\overset{②}{\leq} -\frac{g^\top(\mathbf{H}+\lambda\mathbf{I})^{-1}g}{2} - \left(\frac{2\lambda - 3L\tilde{\varepsilon}}{L}\right)^2 \left(\frac{\lambda}{6} + \frac{L\tilde{\varepsilon}}{3}\right) + \frac{1}{1000\kappa^3 L^2} \\
&\overset{③}{\leq} -\frac{g^\top(\mathbf{H}+\lambda\mathbf{I})^{-1}g}{2} - \frac{2\lambda^3}{3L^2} + \frac{1}{500\kappa^3 L^2} \\
&\leq -\frac{g^\top(\mathbf{H}+\lambda\mathbf{I})^{-1}g}{2} + \frac{1}{500\kappa^3 L^2}
\end{aligned} \tag{D.9}$$
Above, ① is due to Lemma D.2; ② uses our condition on $v$ which gives $L\|v\| \in [2\lambda - 3L\tilde{\varepsilon}, 2\lambda + 3L\tilde{\varepsilon}]$; ③ uses our condition (6.1) on $\tilde{\varepsilon}$.

We now bound $S_1$. For this purpose first we note that if $\lambda_i + \lambda^* \geq \frac{1}{\kappa}$ and $\lambda - \lambda^* \leq \frac{1}{\kappa}$ then
$$2(\lambda_i + \lambda^*) \geq 1/\kappa + \lambda_i + \lambda^* \geq \lambda_i + \lambda \ .$$
Therefore, the sum $S_1$ satisfies
$$S_1 = -\sum_{i:\lambda_i+\lambda^*\geq \frac{1}{\kappa}} \frac{g_i^2}{\lambda_i + \lambda^*} \geq -2\sum_{i:0<\lambda_i+\lambda^*\leq \frac{1}{\kappa}} \frac{g_i^2}{(\lambda_i+\lambda)} \geq -2(g^\top(\mathbf{H}+\lambda\mathbf{I})^{-1}g) \geq 4m(v) - \frac{1}{250\kappa^3 L^2}$$
(Note that we have $\mathbf{H} + \lambda\mathbf{I} \succ 0$.) This finishes the proof of Claim D.3. □

**Claim 6.2.** *If $\lambda_{\min}(\mathbf{H}) \leq -\frac{1}{\kappa}$ then $m(h^*) \geq 1500 \min\left\{m(v), m\left(\frac{\lambda v_{\min}}{2L}\right)\right\} - \frac{1}{500\kappa^3 L^2}$*

*Proof of Claim 6.2.* We derive that
$$\begin{aligned}
m(h^*) &\overset{①}{=} \frac{1}{2}(S_1 + S_2) - \frac{2(\lambda^*)^3}{3L^2} \overset{②}{\geq} \frac{1}{2}(S_1 + S_2) - \frac{16|\lambda_{\min}|^3}{3L^2} \\
&\overset{③}{\geq} 2m(v) - \frac{1}{500\kappa^3 L^2} + 500 \cdot m\left(\frac{\lambda v_{\min}}{2L}\right) - \frac{16|\lambda_{\min}|^3}{3L^2} \\
&\overset{④}{\geq} 2m(v) - \frac{1}{500\kappa^3 L^2} + 1500 \cdot m\left(\frac{\lambda v_{\min}}{2L}\right)
\end{aligned}$$

$$\geq 1500 \min\left\{m(v), m\left(\frac{\lambda v_{\min}}{2L}\right)\right\} - \frac{1}{500\kappa^3 L^2}$$

Above, ① uses equation (D.4), inequality ② follows because we have $\lambda_{\min}(\mathbf{H}) \leq -\frac{1}{\kappa}$ in the assumption and have $\lambda^* \leq -\lambda_{\min}(\mathbf{H}) + \frac{1}{\kappa}$ in the assumption of Case 2 of Main Lemma 1; inequality ③ uses Claim D.3 and Claim D.1; and inequality ④ uses (D.8). This finishes the proof of Claim 6.2. □

**Claim 6.3.** *If $\lambda_{\min}(\mathbf{H}) \geq -\frac{1}{\kappa}$ then $m(h^*) \geq 2m(v) - \frac{16}{\kappa^3 L^2}$*

*Proof of Claim 6.3.* This time we lower bound $S_2$ slightly differently:

$$S_2 \stackrel{①}{\geq} -\frac{4}{\kappa L^2}(\lambda^*)^2 \stackrel{②}{\geq} -\frac{16}{\kappa^3 L^2} \quad (D.10)$$

where ① comes from the second to last inequality from (D.7) and ② comes from $\lambda^* \leq \lambda \leq -\lambda_{\min}(\mathbf{H}) + \frac{1}{\kappa} \leq \frac{2}{\kappa}$ using our assumption in Case 2 of Main Lemma 1.

Putting these together we get that

$$m(h^*) \stackrel{①}{=} \frac{1}{2}(S_1 + S_2) - \frac{2(\lambda^*)^3}{3L^2} \stackrel{②}{\geq} 2m(v) - \frac{1}{500\kappa^3 L^2} - \frac{15}{\kappa^3 L^2} \geq 2m(v) - \frac{16}{\kappa^3 L^2} \ .$$

Above, ① comes from (D.4), ② uses Claim D.3, lower bound (D.10) and $\frac{2(\lambda^*)^3}{3L^2} \leq \frac{16}{3\kappa^3 L^2}$ □

# E Proof of Main Lemma 2

**Main Lemma 2.** *In the same setting as Main Lemma 1, suppose $m(h^*) \geq -\frac{\varepsilon^{3/2}}{300\sqrt{L}}$. Then the output vector $v$ satisfies the following conditions:*

$$\|v\| \leq \|h^*\| + \frac{3}{\kappa L} \quad \text{and} \quad \|\nabla m(v)\| \leq \frac{\varepsilon}{4} + \frac{15}{\kappa^2 L} \ .$$

*Proof of Main Lemma 2.* Let's first note that from the value given in Lemma 5.1,

$$(\lambda^*)^3 \leq \frac{3L^2|m(h^*)|}{2} \leq \frac{L^{3/2}\varepsilon^{3/2}}{200} \ . \quad (E.1)$$

If Case 1 occurs, we have

$$\|v\| \stackrel{①}{\leq} \|(\mathbf{H} + \lambda\mathbf{I})^{-1}g\| + \tilde{\varepsilon} \stackrel{②}{\leq} \frac{2\lambda + 2L\tilde{\varepsilon}}{L} + \tilde{\varepsilon} \stackrel{③}{\leq} \frac{2\lambda^*}{L} + 5\tilde{\varepsilon} \stackrel{④}{\leq} \|h^*\| + \frac{1}{20\kappa L} \ .$$

Above, inequalities ① and ② both use the assumptions of Case 1; inequality ③ uses the fact that $\lambda^* \in [\lambda - L\tilde{\varepsilon}, \lambda + L\tilde{\varepsilon}]$ which again follows from the assumptions of Case 1 (see (D.1)); inequality ④ uses $\|h^*\| = \frac{2\lambda^*}{L}$ from Lemma 5.1 as well as our assumption (6.1) on $\tilde{\varepsilon}$.

As for the quantity $\|\nabla m(v)\|$, we bound it as follows

$$\|\nabla m(v)\| = \left\|g + \mathbf{H}v + \frac{L\|v\|}{2}v\right\| \stackrel{①}{\leq} \|g + (\mathbf{H} + \lambda\mathbf{I})v\| + \lambda\|v\| + L\|v\|^2$$

$$\stackrel{②}{\leq} \|\mathbf{H} + \lambda\mathbf{I}\|\tilde{\varepsilon} + \lambda\|v\| + L\|v\|^2 \stackrel{③}{\leq} (L_2 + 2B)\tilde{\varepsilon} + \frac{\lambda(2\lambda + 3L\tilde{\varepsilon}) + (2\lambda + 3L\tilde{\varepsilon})^2}{L}$$

$$= (L_2 + 2B)\tilde{\varepsilon} + \frac{6\lambda^2}{L} + 15\tilde{\varepsilon}\lambda + 9L\tilde{\varepsilon}^2 \stackrel{④}{\leq} \frac{6(\lambda^* + L\tilde{\varepsilon})^2}{L} + (L_2 + 32B)\tilde{\varepsilon} + 9L\tilde{\varepsilon}^2$$

$$\stackrel{⑤}{\leq} \frac{6(\lambda^*)^2}{L} + (L_2 + 56B)\tilde{\varepsilon} + 15L\tilde{\varepsilon}^2 \stackrel{⑥}{\leq} \frac{\varepsilon}{4} + \frac{15}{\kappa^2 L} \ .$$

Above, inequality ① uses triangle inequality; inequality ② uses $\|v + (\mathbf{H} + \lambda\mathbf{I})^{-1}g\| \leq \tilde{\varepsilon}$; inequality ③ uses $\|\mathbf{H} + \lambda\mathbf{I}\| \leq L_2 + 2B$ and $L\|v\| \leq 2\lambda + 3L\tilde{\varepsilon}$ which comes from our upper bound on $\|v\|$ above; ④ uses the fact that $\lambda^* \in [\lambda - L\tilde{\varepsilon}, \lambda + L\tilde{\varepsilon}]$ which again follows from the assumptions of Case 1 (see

(D.1)); inequality ⑤ uses $\lambda^* \leq 2B$; and inequality ⑥ uses (E.1) together with our assumption (6.1) on $\tilde{\varepsilon}$.

If Case 2 occurs, we have

$$\|v\| \stackrel{①}{\leq} \|(\mathbf{H} + \lambda \mathbf{I})^{-1} g\| + \tilde{\varepsilon} \stackrel{②}{\leq} \frac{2\lambda}{L} + \tilde{\varepsilon} \stackrel{③}{\leq} \frac{2(\lambda^* + 1/\kappa)}{L} + \tilde{\varepsilon} \stackrel{④}{\leq} \|h^*\| + \frac{3}{\kappa L} \quad . \tag{E.2}$$

Above, inequalities ① and ② both use the assumptions of Case 2; inequality ③ uses $\lambda \leq -\lambda_{\min}(\mathbf{H}) + 1/\kappa$ from our assumption of Case 2 as well as $-\lambda_{\min}(\mathbf{H}) \leq \lambda^*$ which comes from Lemma 5.1; inequality ④ uses $\|h^*\| = \frac{2\lambda^*}{L}$ from Lemma 5.1 as well as our assumption (6.1) on $\tilde{\varepsilon}$.

The quantity $\|\nabla m(v)\|$ can be bounded in an analogous manner as Case 1:

$$\|\nabla m(v)\| \leq \|\mathbf{H} + \lambda \mathbf{I}\| \tilde{\varepsilon} + \lambda \|v\| + L\|v\|^2 \leq (L_2 + 2B)\tilde{\varepsilon} + \frac{\lambda(2\lambda + L\tilde{\varepsilon}) + (2\lambda + L\tilde{\varepsilon})^2}{L}$$

$$\stackrel{①}{\leq} \frac{6\lambda^2}{L} + \frac{1}{10\kappa^2 L} \stackrel{②}{\leq} \frac{6(\lambda^* + \frac{1}{\kappa})^2}{L} + \frac{1}{10\kappa^2 L} \leq \frac{12(\lambda^*)^2}{L} + \frac{15}{\kappa^2 L} \stackrel{③}{\leq} \frac{\varepsilon}{4} + \frac{15}{\kappa^2 L} \quad .$$

Above, inequality ① uses our assumption (6.1) on $\tilde{\varepsilon}$; inequality ② uses $\lambda \leq \lambda^* + \frac{1}{\kappa}$ which appeared in (E.2); inequality ③ uses (E.1). $\square$

## F  Proof of Lemma 7.2

**Lemma 7.2.** *The following statements hold for all $i$ until* FastCubicMin *terminates*

(a) $\lambda_i \in [0, 2B]$, $\lambda_i + \lambda_{\max}(\mathbf{H}) \leq 3B$

(b) $\lambda_i + \lambda_{\min}(\mathbf{H}) \geq \frac{3}{10\kappa}$

(c) $\lambda_{i+1} + \lambda_{\min}(\mathbf{H}) \leq \frac{3}{4}(\lambda_i + \lambda_{\min}(\mathbf{H}))$ *unless $\lambda_{i+1} = 0$*

*Moreover when* FastCubicMin *terminates at Line 20 we have $\lambda_i + \lambda_{\min}(\mathbf{H}) \leq \frac{1}{\kappa}$.*

*Proof of Lemma 7.2.* The lemma follows via induction.

To see (a) and (b) at the base case $i = 0$, recall that the definitions of $B$ and $L_2$ together ensure $\lambda_0 + \lambda_{\max}(\mathbf{H}) \leq 3B$ and $\lambda_0 + \lambda_{\min}(\mathbf{H}) \geq \frac{3}{10\kappa}$. Also $\lambda_0 \in [0, 2B]$.

Suppose now for some $i \geq 0$ properties (a) and (b) hold. It is easy to check that $\lambda_i \leq \lambda_{i-1}$ and thus we have $\lambda_i + \lambda_{\max}(\mathbf{H}) \leq 2B$ and $\lambda_i \leq 2B$. This implies property (a) at iteration $i + 1$ also hold. We now proceed to show property (c) at iteration $i$ and property (b) at iteration $i + 1$. Recall that the algorithm ensures

$$\frac{9}{10} \lambda_{\max}((\mathbf{H} + \lambda_i \mathbf{I})^{-1}) \leq w^\top (\mathbf{H} + \lambda_i \mathbf{I})^{-1} w \leq \lambda_{\max}((\mathbf{H} + \lambda_i \mathbf{I})^{-1}) \quad ,$$

and by the definition of $\tilde{w}$ we have

$$\frac{9}{10} \lambda_{\max}((\mathbf{H} + \lambda_i \mathbf{I})^{-1}) - 2\hat{\varepsilon} \leq \tilde{w}^\top w - \hat{\varepsilon} \leq \lambda_{\max}((\mathbf{H} + \lambda_i \mathbf{I})^{-1}) \quad . \tag{F.1}$$

Now, since $\frac{3}{10\kappa} \leq \lambda_i + \lambda_{\min}(\mathbf{H}) \leq 3B$ from the inductive assumption, it follows from the choice of $\hat{\varepsilon}$ that

$$2\hat{\varepsilon} \leq \frac{1}{30B} \leq \frac{1}{10(\lambda_i + \lambda_{\min}(\mathbf{H}))} = \frac{\lambda_{\max}((\mathbf{H} + \lambda_i \mathbf{I})^{-1})}{10} \quad . \tag{F.2}$$

Plugging Equation (F.2) into Equation (F.1) we get

$$\frac{8}{10} \frac{1}{\lambda_i + \lambda_{\min}(\mathbf{H})} = \frac{8}{10} \lambda_{\max}((\mathbf{H} + \lambda_i \mathbf{I})^{-1}) \leq \tilde{w}^\top w - \hat{\varepsilon} \leq \lambda_{\max}((\mathbf{H} + \lambda_i \mathbf{I})^{-1}) = \frac{1}{\lambda_i + \lambda_{\min}(\mathbf{H})} \quad .$$

Inverting this chain of inequalities, we have

$$\frac{\lambda_i + \lambda_{\min}(\mathbf{H})}{2} \leq \Delta \leq \frac{5(\lambda_i + \lambda_{\min}(\mathbf{H}))}{8} \ . \tag{F.3}$$

From this we derive the following implications:

$$\Delta \leq \frac{1}{2\kappa} \implies (\lambda_i + \lambda_{\min}(\mathbf{H})) \leq \frac{1}{\kappa} \tag{F.4}$$

$$\Delta > \frac{1}{2\kappa} \implies (\lambda_i + \lambda_{\min}(\mathbf{H})) > \frac{4}{5\kappa} \tag{F.5}$$

If Condition (F.4) happens, our algorithm FastCubicMin outputs on Line 20; in such a case (F.4) implies our desired inequality $\lambda_i + \lambda_{\min}(\mathbf{H}) \leq \frac{1}{\kappa}$. If Condition (F.5) happens, our choice $\tilde{\lambda}_{i+1} \leftarrow \lambda_i - \frac{\Delta}{2}$ and Equation (F.3) together imply that

$$\frac{3}{4}(\lambda_i + \lambda_{\min}(\mathbf{H})) \geq \tilde{\lambda}_{i+1} + \lambda_{\min}(\mathbf{H}) \geq \frac{11}{16}(\lambda_i + \lambda_{\min}(\mathbf{H}))$$

Combining this with (F.5) we get that

$$\frac{3}{4}(\lambda_i + \lambda_{\min}(\mathbf{H})) \geq \tilde{\lambda}_{i+1} + \lambda_{\min}(\mathbf{H}) \geq \frac{11}{16}\left(\frac{4}{5\kappa}\right) \geq \frac{3}{10\kappa} \ .$$

Therefore, we conclude that property (c) at iteration $i$ holds and property (b) at iteration $i+1$ hold because $\lambda_{i+1} \geq \tilde{\lambda}_{i+1}$. This finishes the proof of Lemma 7.2. □

## G  Proof of Main Lemma 3: Running Time Half

Having proven the correctness of the algorithm, we now aim to bound the overall running time of FastCubicMin, completing the proof of Main Lemma 3. We prove in Appendix H the following lemma:

**Lemma G.1.** *If $\lambda_2 + \lambda_{\min}(\mathbf{H}) \geq c_1 \in (0,1)$ then BinarySearch ends in $O\big(\log(\frac{(\lambda_1 - \lambda_2)B}{c_1 \cdot L \cdot \tilde{\varepsilon}})\big)$ iterations.*

Since in our FastCubicMin algorithm, it satisfies $\lambda_i \leq 2B$ and $\lambda_i + \lambda_{\min}(\mathbf{H}) \geq \frac{3}{10\kappa}$ (see Lemma 7.2), taken together with our choice of $\tilde{\varepsilon}$ we have:

**Claim G.2.** *Each invocation of BinarySearch ends in $O\big(\log(1/\tilde{\varepsilon})\big)$ iterations.*

**Claim G.3.** *FastCubicMin ends in at most $O(\log(B\kappa))$ outer loops.*

*Proof.* According to Lemma 7.2 we have $\frac{3}{4}(\lambda_{i-1} + \lambda_{\min}(\mathbf{H})) \geq \lambda_i + \lambda_{\min}(\mathbf{H})$ so the quantity $\lambda_i + \lambda_{\min}(\mathbf{H})$ decreases by a constant factor per iteration (except possibly $\lambda_i = 0$ the last outer loop in which case we shall terminate in one more iteration). On one hand, we have began with $\lambda_0 + \lambda_{\min}(\mathbf{H}) \leq 3B$. On the other hand, we always have $\lambda_i + \lambda_{\min}(\mathbf{H}) \geq \frac{3}{10\kappa}$ according to Lemma 7.2. Therefore, the total number of outer loops is at most $O(\log(B\kappa))$. □

### G.1  Matrix Inverse

Since the key component of the running time is the computation of $(\mathbf{H} + \lambda_i \mathbf{I})^{-1} b$ for different vectors $b$ we will first bound the condition number of the matrix $(\mathbf{H} + \lambda_i \mathbf{I})^{-1}$ via the following lemma

**Claim G.4.** *Through out the execution of FastCubicMin and BinarySearch whenever we compute $(\mathbf{H} + \lambda_i \mathbf{I})^{-1} b$ for some vector $b$ it satisfies $\frac{\lambda_i + L_2}{\lambda_i + \lambda_{\min}(\mathbf{H})} \leq 10\kappa L_2$.*

*Proof of Claim G.4.* We first focus on Line 5 and Line 11 of FastCubicMin. There are two cases. If $\lambda_i \geq 2L_2$, then according to $-L_2 \mathbf{I} \preceq \mathbf{H} \preceq L_2 \mathbf{I}$ we can bound $\frac{\lambda_i + L_2}{\lambda_i + \lambda_{\min}(\mathbf{H})} \leq 3$ because the left hand side is the largest when $\lambda_i = 2L_2$. If $\lambda_i < 2L_2$, then by Lemma 7.2 we know $\lambda_i + \lambda_{\min}(\mathbf{H}) \geq \frac{3}{10\kappa}$. This implies $\frac{\lambda_i + L_2}{\lambda_i + \lambda_{\min}(\mathbf{H})} \leq 10\kappa L_2$.

We now focus on Line 3 of BinarySearch. We claim that all values $\lambda_{\mathrm{mid}}$ iterated over BinarySearch also satisfy $\lambda_{\mathrm{mid}} + \lambda_{\min}(\mathbf{H}) \geq \frac{3}{10\kappa}$ (because the values $\lambda_{\mathrm{mid}} \geq \lambda_i$ and $\lambda_i$ satisfies $\lambda_i + \lambda_{\min}(\mathbf{H}) \geq \frac{3}{10\kappa}$ according to Lemma 7.2). Therefore, the same case analysis (with respect to $\lambda_{\mathrm{mid}} \geq 2L_2$ and $\lambda_{\mathrm{mid}} < 2L_2$) also gives $\frac{\lambda_i + L_2}{\lambda_i + \lambda_{\min}(\mathbf{H})} \leq 10\kappa L_2$. $\square$

**Claim G.5.** *Line 5 of* FastCubicMin *and Line 3 of* BinarySearch *runs in time* $\tilde{O}(\mathbb{T}_{\mathsf{inverse}}(\kappa L_2, \tilde{\varepsilon}))$.

*Proof.* Whenever we compute $(\mathbf{H} + \lambda_i \mathbf{I})^{-1}b$ for some vector $v$ it satisfies $\|b\| \leq 1/\tilde{\varepsilon}$; therefore to find $v$ satisfying $\|v + (\mathbf{H} + \lambda_i \mathbf{I})^{-1}b\| \leq \tilde{\varepsilon}$ it suffices to find $\|v + (\mathbf{H} + \lambda_i \mathbf{I})^{-1}b\| \leq \tilde{\varepsilon}^2 \|b\|$. This costs a total running time $\tilde{O}(\mathbb{T}_{\mathsf{inverse}}(\kappa L_2, \tilde{\varepsilon}))$ according to Theorem 2.4. $\square$

Therefore by Theorem 2.4, every time we need to multiply a vector $v$ to $(\mathbf{H} + \lambda \mathbf{I})^{-1}$ to error $\delta$, the time required to approximately solve such a system is $\mathbb{T}_{\mathsf{inverse}}(O(\kappa L_2), \delta)$. We will state our running time with respect to $\mathbb{T}_{\mathsf{inverse}}$ as it is the dominant operation in the algorithm.

### G.2 Power Method

We now bound the running time of Power Method in Line 11 of FastCubicMin. It is a folklore (cf. [3, Appendix A]) that getting any constant multiplicative approximation to the leading eigenvector of any PSD matrix $\mathbf{M} \in \mathbb{R}^{d \times d}$ requires only $O(\log d)$ iterations, each computing $\mathbf{M}b$ for some vector $b$. In our case, we have $\mathbf{M} = (\mathbf{H} + \lambda_i \mathbf{I})^{-1}$ so we cannot compute $\mathbf{M}b$ exactly. Fortunately, folklore results on inaccurate power method suggests that, as long as each $\mathbf{M}b$ is computed to a very good accuracy such as $\tilde{\varepsilon}^{-\Omega(\log d)}$, then we can still get a constant multiplicative approximate leading eigenvector that satisfies Line 11 of FastCubicMin. Ignoring all the details (which are quite standard and can be found for instance in [3, Appendix A]), we claim that

**Claim G.6.** *Line 11 of* FastCubicMin *runs in time* $\tilde{O}\left(\mathbb{T}_{\mathsf{inverse}}\left(\kappa L_2, \varepsilon^{-\Theta(\log(d))}\right)\right) = \tilde{O}\left(\mathbb{T}_{\mathsf{inverse}}(\kappa L_2, \varepsilon)\right)$.

### G.3 Lowest Eigenvector

We will now focus on the running time for the computation of the lowest eigenvector of the Hessian which is required in Line 18. We recall Theorem 2.5 from Section 2 which uses Shift and Invert to compute the largest eigenvalue of a matrix.

Since we are concerned with the lowest eigenvector of $\mathbf{H}$ and by assumption $-L_2 \mathbf{I} \preceq \mathbf{H} \preceq L_2 \mathbf{I}$, we can equivalently compute the largest eigenvector of $\mathbf{M} \triangleq \mathbf{I} - \frac{\mathbf{H} + L_2 \mathbf{I}}{2L_2}$ which satisfies $0 \preceq \mathbf{M} \preceq \mathbf{I}$. Note that computing $\mathbf{M}v$ is of the same time complexity as computing $\mathbf{H}v$. By setting $\varepsilon = \delta_\times = \frac{0.01}{\kappa L_2}$ in Theorem 2.5 and running AppxPCA, we obtain a unit vector $w$ such that

$$1 - \frac{w^\top \mathbf{H} w + L_2}{2L_2} = w^\top \mathbf{M} w \geq (1 - 2\delta_\times)\lambda_{\max}(\mathbf{M}) \overset{\text{\textcircled{1}}}{\geq} \lambda_{\max}(\mathbf{M}) - 2\delta_\times \geq 1 - \frac{\lambda_{\min}(\mathbf{H}) + L_2}{2L_2} - 2\delta_\times$$

Above, ① uses $\lambda_{\max}(\mathbf{M}) \leq 1$. Rearranging the terms we obtain $w^\top \mathbf{H} w \leq \lambda_{\min}(\mathbf{H}) + 0.05\kappa$ as desired. In sum,

**Claim G.7.** *The approximate lowest eigenvector computation on Line 18 runs in time* $\tilde{O}\left(\mathbb{T}_{\mathsf{inverse}}(\kappa L_2, \tilde{\varepsilon})\right)$.

### G.4 Putting It All Together

*Running-Time Proof of Main Lemma 3.* Putting together our bounds in Claim G.2 and Claim G.2 which bound the number of iterations, as well as our bounds in Claim G.6, Claim G.5, and Claim G.7 for power method, matrix inverse, and lowest eigenvectors, we conclude that our total running time of FastCubicMin is at most $\tilde{O}(\mathbb{T}_{\text{inverse}}(\kappa L_2, \tilde{\varepsilon}))$, where $\tilde{O}$ contains factors polylogarithmic in $\kappa, L, L_2, B, d$.

By putting together our choice of $\tilde{\varepsilon}$ in Line 2 as well as the running time of either accelerated gradient descent or accelerated SVRG from Theorem 2.4 into formula $O(\kappa L_2, \tilde{\varepsilon})$, we finish the proof of the running time part for Main Lemma 3. □

## H  Proof of Lemma G.1

**Lemma G.1.** *If $\lambda_2 + \lambda_{\min}(\mathbf{H}) \geq c_1 \in (0,1)$ then* BinarySearch *ends in $O\big(\log(\frac{(\lambda_1-\lambda_2)B}{c_1 \cdot L \cdot \tilde{\varepsilon}})\big)$ iterations.*

*Proof of Lemma G.1.* We first note that in all iterations of BinarySearch it always satisfies
$$L\|(\mathbf{H}+\lambda_1\mathbf{I})^{-1}g\| \leq 2\lambda_1 \quad \text{and} \quad L\|(\mathbf{H}+\lambda_2\mathbf{I})^{-1}g\| \geq 2\lambda_2 \ . \tag{H.1}$$
This is true at the beginning. In each of the follow-up iterations, if we have set $\lambda_1 \leftarrow \lambda_{\text{mid}}$ then it must satisfy $L\|v\| + L\tilde{\varepsilon} \leq 2\lambda_{\text{mid}}$ but this implies $L\|(\mathbf{H}+\lambda_{\text{mid}}\mathbf{I})^{-1}g\| \leq 2\lambda_{\text{mid}}$ according to triangle inequality and $\|v+(\mathbf{H}+\lambda_{\text{mid}}\mathbf{I})^{-1}g\| \leq L\tilde{\varepsilon}$; similarly, if we have set $\lambda_2 \leftarrow \lambda_{\text{mid}}$ then it must satisfy $L\|(\mathbf{H}+\lambda_{\text{mid}}\mathbf{I})^{-1}g\| \geq 2\lambda_{\text{mid}}$.

Suppose now the loop has run for at least $\log_2(\frac{\lambda_1-\lambda_2}{\hat{\varepsilon}})$ iterations where $\hat{\varepsilon} \triangleq \frac{L\tilde{\varepsilon}c_1}{40B}$. Then, it must satisfy $\lambda_1 - \lambda_2 \leq \hat{\varepsilon}$. At this point, we compute
$$(\mathbf{H}+\lambda_1\mathbf{I})^{-1} = (\mathbf{H}+\lambda_2\mathbf{I})^{-1} - (\lambda_1-\lambda_2)(\mathbf{H}+\lambda_2\mathbf{I})^{-1}(\mathbf{H}+\lambda_1\mathbf{I})^{-1}$$
and therefore
$$L\|(\mathbf{H}+\lambda_1\mathbf{I})^{-1}g\| \geq L\|(\mathbf{H}+\lambda_2\mathbf{I})^{-1}g\| - L\|(\lambda_1-\lambda_2)(\mathbf{H}+\lambda_2\mathbf{I})^{-1}(\mathbf{H}+\lambda_1\mathbf{I})^{-1}g\|$$
$$\overset{①}{\geq} 2\lambda_2 - \hat{\varepsilon}\|(\mathbf{H}+\lambda_2\mathbf{I})^{-1}\| \cdot 2\lambda_1 \overset{②}{\geq} 2\lambda_1 - 2\hat{\varepsilon} - \hat{\varepsilon}\|(\mathbf{H}+\lambda_2\mathbf{I})^{-1}\| \cdot 2\lambda_1$$
Above, inequality ① uses (H.1) and $\lambda_1 - \lambda_2 \leq \hat{\varepsilon}$; inequality ② uses again $\lambda_1 - \lambda_2 \leq \hat{\varepsilon}$.

Now, we notice that $\|(\mathbf{H}+\lambda_2\mathbf{I})^{-1}\| \leq \frac{1}{c_1}$ and $\lambda_1 \leq 2B$ because $\lambda_2$ only increases and $\lambda_1$ only decreases through the execution of the algorithm. Therefore by the choice of $\hat{\varepsilon} = \frac{\tilde{\varepsilon}c_1}{40B}$, we get
$$L\|(\mathbf{H}+\lambda_1\mathbf{I})^{-1}g\| \geq 2\lambda_1 - L\tilde{\varepsilon}/5 \ .$$
A completely analogous argument also shows that
$$L\|(\mathbf{H}+\lambda_2\mathbf{I})^{-1}g\| \leq 2\lambda_2 + L\tilde{\varepsilon}/5 \ .$$
Therefore, in the immediate next iteration when picking $\lambda_{\text{mid}} \leftarrow (\lambda_1+\lambda_2)/2$, it must satisfy
$$2\lambda_{\text{mid}} - L\tilde{\varepsilon}/2 \leq 2\lambda - L\tilde{\varepsilon}/5 \leq L\|(\mathbf{H}+\lambda_{\text{mid}}\mathbf{I})^{-1}g\| \leq 2\lambda_2 + L\tilde{\varepsilon}/5 \leq 2\lambda_{\text{mid}} + L\tilde{\varepsilon}/2 \ .$$
Then, at this iteration when $v$ is computed to satisfy $\|v+(\mathbf{H}+\lambda_{\text{mid}}\mathbf{I})^{-1}g\| \leq \tilde{\varepsilon}/2$, we also have
$$2\lambda_{\text{mid}} - L\tilde{\varepsilon} \leq L\|v\| \leq 2\lambda_{\text{mid}} + L\tilde{\varepsilon}$$
which means BinarySearch will stop in this iteration. In sum, we have concluded that there will be no more than $O\big(\log(\frac{(\lambda_1-\lambda_2)B}{c_1 \cdot L \cdot \tilde{\varepsilon}})\big)$ iterations. □

## Acknowledgements

We thank Ben Recht for helpful suggestions and corrections to a previous version.